\documentclass[a4paper]{amsart}
\usepackage{amsmath}
\usepackage{rotating}
\usepackage{graphics}

\def\s3{\mathcal{S}_3}
\def\dimo{\textit{Proof}$\quad$}

\def\ss{*}

\def\Torus{\mathbb{T}}

\def\me{\mathbf{g}}
\def\contr{\neg}

\def\ns{\overline{*}}
\def\contr{\rightharpoonup}

\def\R{\mathbb{R}}
\def\C{\mathbb{C}}
\def\Z{\mathbb{Z}}

\def\lefs{\mathit{Lef}^s}
\def\gg{{\mathcal{L}^{s}_{\C}}}
\def\ggp{{\mathcal{L}^{s}_{\C,p}}}
\def\ggr{{\mathcal{L}^s_{\R}}}
\def\ggrp{{\mathcal{L}^s_{\R,p}}}
\def\ggru{{\mathbf{u}\mathcal{L}^s}}
\def\ggrs{{\mathbf{s}\mathcal{L}^s}}
\def\ggrsp{{\mathbf{s}\mathcal{L}^s_p}}
\def\so2{\mathbf{so}(2,\R)}
\def\sl2c{\mathbf{sl}(2,\C)}

\def\ka{K\"{a}hler\;}

\def\bi{{\,\overline{i}\,}}
\def\bj{{\,\overline{j}\,}}
\def\bk{{\,\overline{k}\,}}

\newtheorem{teo}{Theorem}[section]
\newtheorem{cor}[teo]{Corollary}

\newtheorem{lem}[teo]{Lemma}
\newtheorem{pro}[teo]{Proposition}
\newtheorem{dfn}[teo]{Definition}

\newtheorem{rmk}[teo]{Remark}

\title[Natural Lie Algebra bundles on rank two s-K\"ahler manifolds]{Natural Lie Algebra bundles on rank two s-K\"ahler manifolds, abelian varieties and moduli of  curves}
\author{Giovanni Gaiffi, ~Michele Grassi}
\date{}
\begin{document}

\maketitle

\begin{abstract}
We prove that one can obtain natural bundles of  Lie algebras on rank two $s$-\ka manifolds, whose fibres are isomorphic to $\mathbf{so}(s+1,s+1)$, $\mathbf{su}(s+1,s+1)$ and $\mathbf{sl}(2s + 2,\R)$. In the most rigid case (which includes complex tori and abelian varieties) these bundles have natural flat connections, whose flat  global sections act naturally on cohomology.
We also  present several natural examples of manifolds which can be equipped with an $s$-\ka structure with various levels of rigidity: complex tori and abelian varieties, cotangent bundles of smooth manifolds and moduli of pointed elliptic curves.

\end{abstract}
\section{Introduction}

In this paper we  prove that one can obtain natural bundles of Lie algebras on rank two $s$-\ka manifolds, with fibres isomorphic to $\mathbf{so}(s+1,s+1)$ and to  $\mathbf{su}(s+1,s+1)$, $\mathbf{sl}(2s + 2,\R)$. In the most rigid case (which includes complex tori and abelian varieties) these bundles have natural flat connections, whose flat  global sections act naturally on cohomology. \\
An $s$-\ka structure is a direct generalization (with $s$  distinct "\ka forms") of the notion of \ka structure, to which it reduces when $s = 1$.
The original motivation for the introduction in \cite{G1} of $s$-\ka manifolds  was the geometric study of the analytical theory of maps from (open subsets of)  $\R^s$ to a given manifold. Then it was realized that this theory in the case $s = 2$ is well suited for the study of Mirror Symmetry (see \cite{G2},\cite{G3}), as it should be, given that for $s=2$ we are considering maps from Riemann surfaces into general manifolds. In a further specialization of the general features of the theory, already in \cite{G1} and \cite{G2}, and more in detail in \cite{GG1}, it was  argued that, independently of the geometric motivations, one can observe a rich algebraic structure "living" on natural bundles on $s$-\ka manifolds. We decided therefore to embark in a systematic study of these algebraic structures, because of their intrinsic interest, and with the strong belief that if we could master them well enough, we would then be able to apply this theory to interesting problems in algebra and in geometry.\\
The first result of this line of work is the paper \cite{GG2}, where we found  a natural Lie superalgebra bundle on rank three 2-\ka manifolds. We believe that these computations will be very relevant to the geometric study of Mirror Symmetry, and to the search of a geometrical interpolation between the various string theories (see \cite{G2} for a more detailed introduction on this aspect). More specifically, in \cite{G2} it was conjectured that the natural bundles of Lie (super) algebras and of their representations on certain 2-\ka manifolds could provide the natural background on which to build Field Theories; these are rich from the representation theoretic point of view and, once quantized using the language of \cite{G1}, were conjectured to be the right playing field for the search of an M-theory (see for example \cite{DOPW} for a similar approach to the Standard Model in particle physics). A direction more in line with this algebraic study (but actually strongly related to the previous one) is the relationship with Higgs bundles and Hitchin systems (see for example \cite{HT}). \\
In the  present paper  we obtain results which are directly applicable to  questions in Complex geometry and in Algebraic geometry.  The basic reason for this is that a rank two  $s$-\ka  manifold (or a naturally defined double cover of it, in some cases) has a canonical complex structure, with which it becomes \ka of complex dimension $s + 1$. We prove that on these \ka manifolds originating from $s$-\ka geometry there are natural bundles of unitary Lie algebras of signature $(s+1,s+1)$, which act in various ways on differential forms. It is this natural way of representing "large" and well known unitary Lie algebras on differential forms which opens a wide range of geometric applications. For comparison, one should recall that the corresponding constructions for  plain \ka manifolds produce the "Lefschetz" action of $\mathbf{sl}(2,\C)$ on forms and on cohomology, which has a lot of geometric applications and consequences.\\
Let us now introduce  more in detail the geometric and algebraic characters which will play a role. The initial object comes from a generalization and an unconventional point of view on the notion of jet space (see \cite{G1} for details):
\begin{dfn}[\cite{G1}, Definition 2.1 and Corollary 2.6]  A {\em polysymplectic} manifold of rank $r$ is a smooth manifold of dimension $(s+1)r+c$, together with $s$ smooth closed two-forms $\omega_1,...,\omega_s$ such that for any $p\in X$ there is a (Darboux) coordinate  system around $p$ of the form
\[x_1,\ldots,x_r,y^1_1,\ldots,y^1_r,\ldots\ldots,y^s_1,\ldots,y^s_r,z_1,...,z_c\]
for which the forms have the canonical (local) expression
\[\omega_j = \sum_{i=1}^r dx_i\wedge dy^j_i\]
\end{dfn}
\noindent
Of course, when $s = 1$ we recover the usual notion of (possibly degenerate) symplectic manifold. When one adds a Riemannian metric, and asks for the natural compatibility conditions with the polisymplectic data, one comes to our main object of study:\\
\begin{dfn}[\cite{G1}, Definition 7.2]
\label{dfn-ska}
A smooth manifold $M$ of dimension $r(s+1)$ together with a
Riemannian metric $\mathbf{g}$ and $2$-forms
$\omega_{1},...,\omega_{s}$ is {\em $s$-\ka} (of rank $r$) if
for each point of $M$ there exist an open neighborhood $\mathcal{U}$
of $p$ and a system of coordinates
$x_{i},y^{j}_{i}$,$i~=~1,...,r$, $j~=~1,...,s$
on $\mathcal{U}$ such that:\\
1) $\forall j ~ ~\omega_{j}~=~\sum_{i=1}^r dx_{i}\wedge dy^{j}_{i}$,\\
2) $\mathbf{g}_{(\mathbf{x},\mathbf{y})}~=~\sum_{i=1}^r dx_{i}\otimes
dx_{i}~+~\sum_{i,j}dy^{j}_{i}\otimes d
y^{j}_{i}~+~\mathbf{O}(2)$.\\
Any such system of coordinates is called {\em standard}($s$-\ka).
\end{dfn}
\noindent
For $s = 1$ one recovers the usual notion of \ka manifold. As in the case of \ka manifolds, one can use the differential forms associated to the structure to build "wedge" operators on forms, and, using their adjoints, one gets natural Lie algebras. Again as in the case of \ka manifolds, to build correctly a theory involving also the adjoints of such wedge operators, it is necessary to consider their pointwise action, and to recover the global operators on forms as global sections of corresponding bundles of Lie algebras. When $s= 1$ one obtains the classical $\mathbf{sl}(2,\C)$ action on the forms of a \ka manifold (and on its cohomology using the Hodge identities). In the case $s > 1$  there is a qualitatively different situation, in that there are more natural differential two-forms than one could initially guess. Indeed, in addition to the structural forms $\omega_1,...,\omega_s$ which generalize directly the \ka form, there are also "mixed" forms $\omega_{jk}$ for any pair of indices $j,k\in\{0,...,s\}$, including the structural ones via the identifications
\[\omega_j ~=~ \omega_{0j}~\quad\text{for}~j\in\{1,...,s\}\]
The precise description of these derived natural forms will be given in the next section. Here however we can already use them to build corresponding "wedge" operators, Lefschetz style:
\begin{dfn}
\label{dfn:Loperators}
For $\phi\in
\Omega^*_\C X$ and $j,k\in\{0,\ldots s\}$ with $j\not= k$,
\[L_{jk}(\phi) ~=~ \omega_{jk}\wedge \phi  ~=~ -L_{kj}(\phi)\]
\end{dfn}
Some canonical mutually orthogonal distributions \(W_i\)  (\(i=0, 1,2,\ldots , s\)) are induced on $T_p^* X$ by the forms \(\omega_{jk}\) (see Section \ref{sec:defi}).
Therefore other  natural operators, called \(V_i\)  (\(i=0, 1,2,\ldots , s\)) arise from wedging with the  local volume forms  of these distributions.

One then uses all  these operators, and their (pointwise) adjoints,  to build a  natural bundle of real  Lie algebras \(\ggr\) (and its complexified bundle \(\gg\))   acting on forms. To be precise, one can fix a point $p\in X$, and on this point one can restrict the actions above, to obtain bundles of Lie algebras on the $s$-\ka manifold. This approach has many advantages, among which the fact that these bundles will exist also in situations in which the single operators used to define them do not have global sections on all of $X$. \\

Coming to a more detailed description of the contents of the present paper, in Section \ref{sec:defi}  we give the definiton of {\em almost} $s-$\ka structure, which is a weaker version of the   defintion of $s-$\ka structure. Then we  provide  a first geometric description of an (almost) $s-$\ka  manifold \(X\): we  discuss the existence of an (almost) complex structure, the natural distributions on the  cotangent space, the group of local structure-preserving tranformations and  the orientability  properties.

Section \ref{sec:natural} is devoted to the definition of  the  bundles of    Lie algebras \(\ggr\) and \(\gg\) on \(X\)  which constitute the main object to be studied in this paper.   We also  point out two  other real forms (\( \ggrs, \ggru\)), defined in terms of geometric generators,  of the  bundle \(\gg\).  We then   define  \(\lefs\) as the  real sub-bundle of \(\gg\) which is the direct generalization of  the classical $\mathbf{sl}(2,\C)$ Lefschetz bundle of  \ka geometry.

The sections from 4 to 7   are a detailed study of the fibres of  above mentioned Lie algebra bundles: in Section 4
the Lefschetz bundle \(\lefs\)   is studied  in detail, by showing   some fundamental relations among its generators;  the fibres of the bundle turn out to be   isomorphic to the orthogonal algebra \(\mathbf{so}(s+1,s+1,\R)\) and Serre generators are presented in terms of simple brackets of geometric generators (Theorem \ref{pro:dn}).

Sections 5 and 6 are devoted to the complete description of the main complex bundle  \(\gg\):  we use  the Hodge decomposition on \({\bigwedge^*}_\C T_p^* X\)   with respect to the (almost) complex structure and Clifford algebra techniques to show that the fibres of \(\gg\) are isomorphic to \(\mathbf{sl}(2s+2,\C)\); furthermore, we   characterize \(\gg\) as the bundle of all quadratic elements  of trace zero (compatible with the almost complex structure)  of a Clifford algebra bundle (Theorem \ref{teo:slstructure}).

Section 7 focusses on the real forms of \(\gg\) which turn out to be  interesting both from the algebraic and the geometric point of view. In fact \(\ggrs\) is proven to be the bundle  of the real split form of \(\gg\), while \(\ggr\) and \(\ggru\) are shown to be the real bundles of  operators which preserve two natural non degenerate hermitean inner products on \({\bigwedge^*}_\C T_p^* X\) (Theorems \ref{pro:realform}, \ref{teo:ggru}, \ref{teo:ggr}) . A superHermitean variant of one of these inner products  was introduced in \cite{GG2} to study  rank three WSD structures.  A computation of the signature shows that the fibres  of $\ggr$ and $\ggru$ are unitary Lie algebras  isomorphic to \(\mathbf{su}(s+1, s+1)\).  We observe that     the complete description of \(\ggr\) fully answers to the  question (first rised in \cite{G1} and then more precisely in the rank two case in \cite{GG1}) on the nature of the algebraic bundles  generated  by the real  canonical operators associated to an \(s\)-\ka structure.  Furthermore, at the end of the section, \(\lefs\) is shown to coincide with \(\ggr\cap \ggrs\).

After this presentation of the natural  Lie bundles of an \(s\)-\ka structure, we devote  Sections \ref{sec:tori} and \ref{sec:moduli} to the construction of some examples of (full, almost, or pointwise) \(s\)-\ka structures.  In Section \ref{sec:tori} we first recall the standard examples of \cite{G1} built using iterated cotangent bundles of smooth Riemannian manifolds. We then observe how any real torus of dimension $r(s+1)$ can be given many non-equivalent translation invariant $s$-\ka structures. In Corollary \ref{cor:cohomology} we prove that when $X$ is compact orientable $s$-\ka   there is a natural action of the flat global sections of the bundles of Lie algebras  $\gg,\ggr,\ggrs,\ggru,\lefs$ on cohomology. Then, in Theorem \ref{teo:rivfinito}, we characterize the rank two $s$-\ka manifolds which are quotients of tori in terms of the Calabi-Yau condition. In Section \ref{sec:moduli} we put a pointwise rank two  $s$-\ka structure on the moduli space  of elliptic curves with $s+2$ punctures, which depends naturally on any chosen \ka metric.

\section{Introduction to the geometric setting}
\label{sec:defi}
This section is again introductory in nature, but with a stronger emphasis on the geometric aspects of the theory. It would be too long to describe all the general facts on $s$-\ka geometry here, so we will list here only the most relevant ones for our purposes, while referring to \cite{G1} for a more thorough analysis. First, as will have been clear already to the reader, one can isolate the pointwise aspects of the definition of an $s$-\ka manifold. The notion of almost $s$-\ka manifold given below is actually a hybrid between pointwise and local properties, which was introduced in \cite{G1} (in the nondegenerate case) in the belief that this mix could be best suited to the purposes of that paper:

\begin{dfn} 
\label{dfn:aska} A {\em almost $s$-\ka structure of rank $r$} on a  smooth manifold $X$ of dimension $(s+1)r+c$ is given by a Riemannian metric $\mathbf{g}$ and a smooth differential two-form $\omega_{jk}$ for any pair of indices $j,k\in\{0\ldots s\}$ such that  $\omega_{01},...,\omega_{0s}$ give a polysymplectic structure, and for any point $p\in X$ there is an orthonormal basis of $T_p^* X$ made of covectors $v_{ij}$ for $i\in\{1,..,r\}$ and $j\in\{0,...,s\}$ and $u_1,...,u_c$  such that, ,
\[\omega_{jk} = \sum_{i=1}^r v_{ij}\wedge v_{ik}\]
The forms $\omega_{jk}$ with $j,k\in\{1,...,s\}$ are called {\em dualizing forms}.
\end{dfn}
\begin{rmk} One can directly check that giving an almost $s$-\ka structure is equivalent to giving the forms $\omega_{01},...,\omega_{0s}$ and a compatible metric, which is exactly what is needed for an $s$-\ka structure, except for the local condition 2), which is equivalent to  the invariance of the forms with respect to the Levi-Civita connection. This is the way in which $s$-\ka manifolds were introduced in \cite{G1}.
\end{rmk} 
 Clearly there are some redundancies in the definition given above: for example, one has always $\omega_{jk} = - \omega_{kj}$. Observe also that an almost 1-\ka manifold is simply an almost \ka manifold, and for this reason in this paper we consider only the case $s\geq 2$ which is moreover the range where our constructions do exist. Recall also that an almost 2-\ka manifold in which the structure forms are closed is a Weakly Self Dual manifold (see \cite{G2}, Definition 2.6), or WSD manifold for short.\\
{\em 
For the algebraic constructions to be discussed in this paper, all that is needed is a rank 2 almost $s$-\ka structure (actually for the main construction we will need only the pointwise part of the definition)}.\\
 The almost $s$-\ka structure on a manifold $X$ splits its cotangent space as
$T_p^*X = W_0 \oplus W_1 \oplus\cdots\oplus W_s$ where the $W_j$ are  $s+1$
mutually orthogonal canonical distributions defined as:
\[W_j = \{\phi\in T_p^*X~|~\phi\wedge \omega_{jk} = 0~~\text{for}~k~\text{in}~0,..,\hat{j},..,s\}\]
The almost $s$-\ka structure also determines canonical pairwise linear
identifications among the  spaces $W_j$, so that one can also
write $T_p^*X \cong W_0\otimes_\R \R^{s+1}$ or more simply
\[T_p^*X \cong W\otimes_\R \R^{s+1}\] where $W = W_0 \cong W_1 \cong \cdots \cong W_s$.\\
Let us now come back to the canonical  operators $L_{jk}$ mentioned in
the Introduction. \\
We now choose an orientaion of $W_0$ at a fixed point $p\in X$, and a (non-canonical) orthonormal basis
$\gamma_1,\gamma_2$ compatible with this orientation;  this together with the standard
identifications of the $W_j$ determines an orientation and an orthonormal basis for
$T_p^*X$, which we write as $\{v_{ij}= \gamma_i\otimes
e_j~|~i=1,2,~j=0,..,s\}$.  We remark that the $v_{ij}$ are an
\textit{adapted coframe} for the almost $s$-\ka structure, and therefore we
have the explicit expressions:
\[\omega_{jk} = v_{1j}\wedge v_{1k} + v_{2j}\wedge v_{2k}\]
A different choice of the $\gamma_1,\gamma_2$ would be related to
the previous one by an element in $\mathbf{O}(2,\R)$. The Lie algebra of the group
$\mathbf{O}(2,\R)$ expressing the change from one adapted
basis to another is generated point by point by the
operator $J$, which is determined and determines a (pointwise, local or global if possible) orientation of the distribution $W_0$:
\begin{dfn}
\label{dfn:so2R} The operator $J\in End_\R(\bigwedge^*T^*_p(X))$ associated to the standard basis $v_{ij}$ is defined as
\[J(v_{1j}) = v_{2j},\qquad J(v_{2j}) = - v_{1j}\qquad \text{for}~
j \in \{0,1,\ldots,s\}\] and $J(v\wedge w) = J(v)\wedge w + v \wedge J(w)$
for $v,w \in \Lambda^*T^*_pX$
\end{dfn}
\begin{rmk}
\label{rmk:j}
As $J$ commutes with itself, and it is determined at every point $p\in X$ by an orientation of $(W_0)_p\subset T_pX$,  it is always well defined locally. Of course, $J$ admits a global determination if and only if $W_0$ admits a global orientation. This happens for example if $X$ is orientable and $s$ is even.
\end{rmk}
Whenever we will need a local volume form on $X$, we will use the one induced by a local choice of $J$ which we will call  $\Omega_p$ over the point $p\in X$.\\
From the above considerations it follows the following fundamental remark:
\begin{rmk} An (almost) $s$-\ka manifold of rank 2 is in particular an (almost) complex manifold of complex dimension s+1, when there is a global determination of $J$. This happens in particular when $X$ is orientable and $s$ is even.
\end{rmk}
For this reason, rank two (almost) 2-\ka manifolds can be seen as a chapter in (almost) complex geometry. This allows on one hand to "import" the techniques of complex geometry into the realm of almost 2-\ka geometry, and on the other hand allows one to apply the results of almost  2-\ka geometry to the complex world.  Summing up, we have that
 \[\text{$s$-\ka} \implies \text{almost $s$-\ka}\implies \text{Polysymplectic}\]
and in the rank two case we have moreover that locally
 \[\text{(almost) $s$-\ka} \implies \text{(almost) \ka}\]

When the structure is $s$-\ka, one has that all the structure forms are covariant constant with respect to the Levi-Civita connection associated to the metric. This allows one to perform many of the same constructions that one usually performs in the \ka case. In particular, one recovers (the analog of) the Hodge identities, and the adjoints of the canonical operators $L_{jk}$ operate on cohomology (see Theorem \ref{teo:skaid} and Corollary \ref{cor:cohomology}). This is the context in the case of Abelian varieties, which in our opinion will provide many interesting applications of the constructions to be detailed in the present paper.

For a general  rank of the structure  $r\geq 1$, many of the above considerations generalize; for example  the group of pointwise transformations which preserve the structure is $\mathbf{O}(r)$. As we have seen above, in the $r = 2$ case we obtain $\mathbf{O}(2)$ whose algebra is generated  by $J$, while the $r = 3$ case (in which comes into play $\mathbf{O}(3)$) was discussed in detail in \cite{GG2}. Clearly however, not everything generalizes  to arbitrary rank: for example, a rank three $s$-\ka manifold may be of (real) dimension $9$, which is odd and therefore it is impossible to have an almost complex structure on such a manifold. Still in case $r = 3$, one has natural operators also in odd degree, and therefore the natural algebras which come out of the geometry are Lie superalgebras, instead of Lie algebras (see \cite{GG2}).

\section{Construction of the natural algebras}
\label{sec:natural}
In this section we fix a point $p$ in an almost $s$-\ka manifold $X$ and we mostly work on tensor powers of $T_p X$.

As was mentioned in the previous sections, using the forms $\omega_{jk}$ of the almost $s$-\ka structure, we can build corresponding operators on forms, much in the way as the $L$ operator is built on \ka manifolds:\\

\noindent
\textbf{Definition~\ref{dfn:Loperators}}\textit{ For $\phi\in
\Omega^*_\C X$ and $j,k\in\{0,\ldots s\}$ with $j\not= k$,
\[L_{jk}(\phi) ~=~ \omega_{jk}\wedge \phi  ~=~ -L_{kj}(\phi)\] }\\
The above operators restrict also to $\bigwedge^*T^*_pX$ for any $p\in X$ where, using the chosen (orthonormal) basis, one can define also corresponding (non canonical) wedge and contraction operators:
\begin{dfn} Let $i\in\{1,2\}$, $j\in\{0,1,\ldots,s\}$ and $p\in X$. The operators $E_{ij}$ and $I_{ij}$ are respectively the wedge and the contraction operator with the form $v_{ij}$ on $\bigwedge^*T^*_pX$ (defined using the given basis); we use the notation $\frac{\partial}{\partial v_{ij}}$ to indicate the element of $T_pX$ dual to $v_{ij}\in T^*_pX$:
\[E_{ij}(\phi) = v_{ij}\wedge \phi,\qquad I_{ij}(\phi) = \frac{\partial}{\partial v_{ij}}\contr
\phi\]
\end{dfn}
\begin{pro}
\label{pro:clifrelations}
The operators $E_{ij},I_{ij}$ satisfy the following relations:
\[ \forall i,j,k,l\qquad E_{ij}E_{kl} = -E_{kl}E_{ij},\quad I_{ij}I_{kl} = - I_{kl}I_{ij}\]
\[\forall i,j \qquad E_{ij}I_{ij} + I_{ij}E_{ij} = Id\]
\[\forall (i,j)\not= (k,l) \qquad E_{ij}I_{kl} = - I_{kl}E_{ij}\]
\[\forall i,j\qquad E_{ij}^* = I_{ij},\quad I_{ij}^* = E_{ij}\]
where $*$ is adjunction with respect to the metric.
\end{pro}
\dimo The proof is a simple direct verification, which we omit. \qed\\
It is then immediate to check that:
\begin{pro}
\label{pro:dfnJ} $J$  can be expressed on the whole $\bigwedge^*T^*_pX$ as
\[J = \sum_{j=0}^s\left(E_{2j}I_{1j} - E_{1j}I_{2j}\right)\]
\end{pro}
\begin{rmk}
\label{rmk:dfnJ}
From this expression and the previous proposition one obtains that  $J^* = -J$, i.e.  for every $p$ the Lie algebra generated
by $J$ is a subalgebra of $\mathbf{o}(\bigwedge^*T^*_pX)$ isomorphic to $\so2\cong \R$. 
\end{rmk}
Using the (non canonical) operators $E_{ij}$ we can obtain simple expressions for the
pointwise action of the canonical wedge operators $V_j$ associated to the volume forms of the distributions $W_j$:
\begin{dfn}
\label{dfn:Voperators} For $\phi\in \bigwedge^*T^*_pX$ and $j\in\{0,\ldots,s\}$,
\[V_j(\phi) = E_{1j} E_{2j}(\phi)\]
\end{dfn}
Remember however that the operators $V_j$, being  simply multiplication by the volume forms of the spaces $W_j$,  depend on the choice of a pointwise orientation for these spaces, which is implied for example by the choice of a determination for the operator $J$. Notice that when $s$ is even, and $X$ is oriented, it is always possible to define $J$ (and consequently $V_j$) globally on $X$. On the opposite extreme situation, if $X$ is non-orientable, it is certainly not possible to orient globally any one of the distributions $W_j$ (and a fortiori you cannot determine $J$ globally).

The riemannian metric induces a Riemannian metric on $T^*_p X$ and
on the space $\bigwedge^*T^*_pX$.
\begin{dfn}
\label{dfn:lambda-A-operators} For $j\not=k \in\{0,1,\ldots,s\}$
\[\Lambda_{jk} = L_{jk}^*,\qquad A_j =
V_j^*\]
\end{dfn}
By
construction the canonical operators $L_{jk},\Lambda_{jk}$ on
$\bigwedge^*T^*_pX$ are the pointwise restrictions of
corresponding global operators on smooth differential forms, which
we indicate with the same symbols: for $j\not= k \in \{0,1,\ldots,s\}$,
\[L_{jk},\Lambda_{jk}: \Omega^*(X)\to
\Omega^*(X)\]
In the study of \ka geometry, a central role is played by the Lie algebra generated by Lefschetz operator and its adjoint. The direct generalization of that algebra to the setting of (almost) $s$-\ka manifolds is the following:
\begin{dfn}
\label{dfn:lefs} The smooth bundle of Lie algebras $\lefs$ is the real sub-bundle of Lie
algebras of $End_\R\left(\Omega^*(X)\right)$ generated locally by the
operators
\[\{L_{jk},\Lambda_{jk}~|~\text{for}~ ~ j = 0,1,\ldots,s\}\]
\end{dfn}
The $V_j,A_j$ instead can be always determined locally via a local determination of the operator $J$ even when $s$ is odd. 
Summing up:
\begin{dfn}
\label{dfn:algebragg} The smooth bundle of Lie algebras $\ggr$ is the real sub-bundle of Lie
algebras of $End_\R\left(\Omega^*(X)\right)$ generated locally by the
operators
\[\{L_{jk},V_j,\Lambda_{jk},A_j~|~\text{for}~ ~ j = 0,1,\ldots,s\}\]
for any fixed determination of $J$.
The $*$-Lie algebra $\gg$ is $\ggr\otimes_\R \C$.
The $*$ operator on $\gg$ is induced by the adjoint with respect to
the Hermitean metric induced by the Riemannian one via complexification. 
\end{dfn}
As mentioned in the Introduction, in the present paper we will describe completely the structure of the fibers of the bundles $\lefs,\ggr,\gg$, and we will further describe two other real forms of $\gg$, which are especially significant from a geometric point of view. Here are their definitions:
\begin{dfn}
\label{dfn:ggrs}
The real form $\ggrs$ of the complex bundle of $*$-Lie algebras $\gg$ is generated (as a bundle of real Lie algebras) by the local operators:
\[L_{jk},~\imath V_j,~\Lambda_{jk},~\imath A_j\]
\end{dfn}
\begin{dfn}
\label{dfn:ggru}
The real form $\ggru$ of the complex bundle of $*$-Lie algebra $\gg$ is generated (as a bundle of real Lie algebras) by the local operators:
\[iL_{jk},~\imath V_j,~i\Lambda_{jk},~\imath A_j\]
\end{dfn}

\section{Clifford algebras and a natural presentation of  $\lefs$ as a \(\mathbf{so}(s+1,s+1,\R)\) bundle}
\label{sec:DN}
In this section we will show that $\ggr$ lies inside a (real) Clifford algebra bundle over the $(4s + 4)$-dimensional real bundle $TX\oplus T^* X$;  we will also point out that the natural bundle of Lie subalgebras  \(\lefs\subset \ggr\) is isomorphic to the constant bundle having as fibre the orthogonal algebras \(\mathbf{so}(s+1,s+1,\R)\). Notice that the above considerations do not apply to the $s=1$ (\ka) situation; $\lefs$ in that case is simply a constant $\mathbf{sl}(2,\R)$ bundle,  as it is well know classically. Notice also that this global trivialization of $\lefs$ does not depend on a determination of the (almost) complex structure $J$. \\
In the following we define some new operators, and in the meantime we introduce a unifying notation which concerns the $L_{jk},\Lambda_{jk}$. These operators will be shown in Corollary ~\ref{cor:ibj} to be (global) sections of $\ggr$. 
\begin{dfn} For $j,k\in\{0,\ldots,s\}$
\[L_{jk}=\sum_{i = 1}^2 E_{ij}E_{ik} \qquad  L_{j\bar{k}}=\sum_{i = 1}^2 E_{ij}I_{ik} \]

\[L_{\bar{k}\bar{j}}=\Lambda_{jk}=\sum_{i = 1}^2 I_{ik}I_{ij} \qquad  L_{\bar{j}k}=\sum_{i = 1}^2 I_{ij}E_{ik}\]
In accordance with the notation introduced in \cite{G2} Section 7, we will use the shortcuts $L_{\alpha\beta}$ with $\alpha,\beta\in\{0,\ldots,s,\overline{0},\ldots,\overline{s}\}$, with the convention that $\overline{\overline{\alpha}} = \alpha$.
\end{dfn}
Notice that with the above notation $L_{\alpha\alpha} = 0$ for any $\alpha\in\{0,\ldots,s,\overline{0},\ldots,\overline{s}\}$.
\begin{lem}
\label{lem:relaz}
Given   \(\alpha,\beta,\gamma \in\{0,\ldots,s,\overline{0},\ldots,\overline{s}\}\) with $\alpha\not=\beta,\alpha\not=\overline{\gamma},\gamma\not=\overline{\beta}$ :

\[[L_{\alpha\beta},L_{\overline{\beta}\gamma}] = L_{\alpha\gamma}\]
Given   \(\alpha\not= \beta \in\{0,\ldots,s,\overline{0},\ldots,\overline{s}\}\):
\[[L_{\alpha\beta},L_{\overline{\beta}\,\overline{\alpha}}] = L_{\alpha\overline{\alpha}} + L_{\beta\overline{\beta}}\]
Given   \(\alpha,\beta,\gamma,\delta \in\{0,\ldots,s,\overline{0},\ldots,\overline{s}\}\) with $\{\alpha,\beta\}\cap\{\overline{\gamma},\overline{\delta}\} = \emptyset$ :

\[[L_{\alpha\beta},L_{\gamma\delta}] = 0\]
\end{lem}

\dimo We prove the first relations with $\alpha = i,\beta = j,\gamma = \overline{k}$ and the second ones with $\alpha = i,\beta = j$. The other cases of the first and second are proved exactly with the same passages. The third set of relations is straightforward due to the anticommutativity of the degree one operators which appear in the expressions of $L_{\alpha\beta},L_{\gamma\delta}$ .\\ 
For the first set of relations, a direct computation  which is based on  the fundamental relations \ref{pro:clifrelations} among the operators \(E_{ij}\) and \(I_{rs}\) proves:
\[ [L_{ij}, L_{\bj \bk}]=\sum_r E_{ri}E_{rj} \sum_s I_{sj}I_{sk}- \sum_s I_{sj}I_{sk} \sum_r E_{ri}E_{rj} =\] \[= \sum_r E_{ri}E_{rj} \sum_s I_{sj}I_{sk}- \sum_{s\neq r} E_{ri}E_{rj}I_{sj}I_{sk} - \sum_{s \; ( s= r)} I_{sj}I_{sk}E_{si}E_{sj}=\]
\[= \sum_r E_{ri}E_{rj} \sum_s I_{sj}I_{sk}- \sum_{s\neq r} E_{ri}E_{rj}I_{sj}I_{sk} + \sum_{s \; ( s= r)} E_{si} I_{sj}E_{sj}I_{sk}=\]
\[= \sum_r E_{ri}E_{rj} \sum_s I_{sj}I_{sk}- \sum_{s\neq r} E_{ri}E_{rj}I_{sj}I_{sk} + \sum_{s } E_{si} I_{sk}-\sum_{s } E_{si} E_{sj}I_{sj}I_{sk}=\]
\[=\sum_{s } E_{si} I_{sk} =L_{i\bk}\]

The second set of relations is proved as follows:
\[ [L_{ij}, L_{\bj \bi}]=\sum_r E_{ri}E_{rj} \sum_s I_{sj}I_{si}- \sum_s I_{sj}I_{si} \sum_r E_{ri}E_{rj} =\] \[= \sum_r E_{ri}E_{rj} \sum_s I_{sj}I_{si}- \sum_{s\neq r} E_{ri}E_{rj}I_{sj}I_{si} - \sum_{s \; ( s= r)} I_{sj}I_{si}E_{si}E_{sj}=\]
\[= \sum_r E_{ri}E_{rj} \sum_s I_{sj}I_{si}- \sum_{s\neq r} E_{ri}E_{rj}I_{sj}I_{si} - \sum_{s } E_{si}E_{sj}I_{sj}I_{si}+\sum_{s} E_{sj}I_{sj}+\sum_{s} E_{si}I_{si}=\]\[=L_{i\bi}+L_{j\bj}\]
\qed\\
\begin{cor}
Given any choice of indices \(j\neq k\), the elements \(L_{\bj k}, L_{j\bk}\) belong to \(\Gamma(X,\lefs)\subset \Gamma(X,\ggr\cap \ggrs)\). Furthermore, for every \(j=0,1,2,\ldots , s\), the elements \(L_{j\bj}\) belong to \(\Gamma(X,\lefs)\subset \Gamma(X,\ggr\cap \ggrs)\).
\label{cor:ibj}

\end{cor}
\dimo For any fixed $p\in X$, the values of the elements  \(L_{jk}\)  and  \(L_{\bj\bk}\) at $p$ are (maybe up to a sign) among the generators of  the fibre of \(\ggr\cap\ggrs\) at $p$.  To show that \(L_{j\bk}\) is a section of  \(\ggr\cap\ggrs\) we notice that, since \(s\geq 2\), we can find an index \(i\in \{0,1,2,\ldots ,s\}\) which is different from both \(j\) and \(k\).  Then we can  use the lemma above and construct \(L_{j\bk}\) as:
\[[L_{ji},L_{\bi\bk}]=L_{j\bk}\]
The element $L_{\bj k}$ is equal to  $-L_{j \bk}^*$ and therefore also is a section of $\ggr\cap\ggrs$.
As for the last assertion, it follows from the first one and the fact that (according to the above lemma)
\([L_{ij}, L_{\bj \bi}]=L_{i\bi}+L_{j\bj}\) and \([L_{i\bj}, L_{j \bi}]=L_{i\bi}-L_{j\bj}\). 
\qed \\

The operators defined below give rise to a set of Serre generators for \(\Gamma(X,\lefs )\), as shown in the following Theorem.

\begin{dfn}
\label{dfn:ei}
Let us define:\\
\(e_1=L_{1\overline{0}}\)\\
\(e_2=L_{2\overline{1}}\)\\
\(e_3=L_{3\overline{2}}\)\\
.....\\
\(e_{s-1}=L_{s-1\overline{s-2}}\)\\
\(e_{s}=L_{s\overline{s-1}}\)\\
\(e_{s+1}=L_{\overline{s-1}\overline{s}}\)\\
Moreover, for every \(i=1,2,\ldots, s+1\),  let \(f_i\) be the adjoint of \(e_i\).
\end{dfn}
\begin{teo}
The global operators \(e_i\), \( f_j\) and \(h_i=[e_i,f_i]\)  restrict to a set of Serre generators of $\lefs_p$ for any $p\in X$, and    $\lefs$  is (canonically) a trivial Lie algebra bundle with fibre  isomorphic to \(\mathbf{so}(s+1,s+1,\R)\).
\label{pro:dn}
\end{teo}

\begin{figure}[h]
  \begin{center}
    \includegraphics{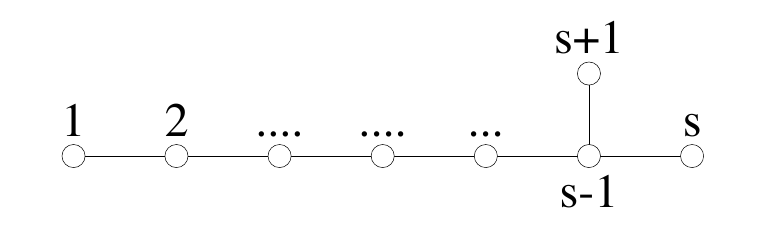}
  \end{center}
  \caption{\label{DynkinD} 
           The Dynkin diagram of  type  $\mathbf{D}_{s+1}$ with labels.}
\end{figure}

\dimo  From the previous corollary, the global operators \(e_i\) , \( f_j\) and \(h_i=[e_i,f_i]\) are sections of $\lefs$. It is immediate, using Lemma \ref{lem:relaz}, to check that these elements are also enough to produce a set of linear generators of $\lefs_p$. We are left with the verification of the Serre relations  for a root system  of type \(\mathbf{D}_{s+1}\). We consider  a  basis of   simple roots   \(\alpha_1,\alpha_2,\ldots ,\alpha_{s-2}, \alpha_{s-1}, \alpha_{s},\alpha_{s+1}\) indexed according to the  labelled Dynkin diagram in 
Figure \ref{DynkinD}, and think of  the operator \(e_i\) (resp. \(f_i\)) as a generator of the root space associated to \(\alpha_i\) (resp. \(-\alpha_i\)).

Then we have to verify that the following relations hold in \(\ggr\):
\begin{enumerate}
\item \([h_i,h_j]= 0\)
\item  \([h_i,e_i]=2e_i\),  \([h_i,f_i]=-2f_i\)
\item \((ad \ e_i) ^{1-\alpha_i(h_j)}e_j =0\) and \((ad \ f_i) ^{1-\alpha_i(h_j)}f_j =0\) for \(i\neq j\).
\item \([e_i,f_j]=0\) for \(i\neq j\).
\item \([h_i, e_j]=\alpha_i(h_j)e_j\),  \([h_i, f_j]=-\alpha_i(h_j)f_j\)

\end{enumerate}
From the relations above if follows that  the $h_1,\ldots,h_{s+1}$ span a Cartan subalgebra of the real Lie algebra generated by the $e_i,f_j$, with real eigenvalues. This proves that the algebra is the split real form $\mathbf{so}(s+1,s+1,\R)$ of $\mathbf{so}(2s +2,\C)$.\\
Concerning the proof of these relations, they are actually all consequence of Lemma \ref{lem:relaz}.
Relations of type \((2)\), for instance,  are all verified using the same computation, which  we  show in the example of    \[[h_1,e_1]=2e_1\] 
This follows from the observation that \(h_1=[L_{1,\overline{0}},L_{0,\overline{1}}]=L_{1,\overline{1}}-L_{0,\overline{0}}\)
and then
\[[h_1, e_1]=[L_{1,\overline{1}}-L_{0,\overline{0}}, L_{1,\overline{0}}]= L_{1,\overline{0}}- [L_{0,\overline{0}}, L_{1,\overline{0}}]=2 L_{1,\overline{0}}=2e_1\]
Among the  last relations to be verified  we show as final examples:
\[[h_{s+1}, e_s]=[[L_{\overline{s-1}\overline{s}},L_{s(s-1)}], L_{s\overline{s-1}}]=\]

\[= - [L_{s-1 \overline{s-1}}+L_{s\overline{s}}, L_{s\overline{s-1}}]= L_{s\overline{s-1}}- L_{s\overline{s-1}}=0\]
and, again by Lemma \ref{lem:relaz}, 
\[[h_{s+1}, e_{s-1}]= - [L_{s-1 \overline{s-1}}+L_{s,\overline{s}},L_{s-1\overline{s-2}}]=- [L_{s-1 \overline{s-1}},L_{s-1\overline{s-2}}]-0=\] 
\[-L_{s-1\overline{s-2}}=-e_{s-1}\]
\qed \\
\begin{rmk} We notice that  Theorem \ref{pro:dn}  is in accordance with   \cite{GG1} and \cite{GG2} where   the specialization of these computations to the case of   \(WSD\) manifolds  of rank two and three    led us to the description of a natural subalgebra isomorphic to  \(\mathbf{sl}(4,\R) \cong \mathbf{so}(2,2,\R)\).
\end{rmk}
An alternative interpretation of the relations in Lemma \ref{lem:relaz} and of the appearance of $\mathbf{D}_{s+1}$ is thruough the use of two different Clifford Algebras, which will play a prominent role in the rest of this paper. For the first one, generalizing to arbitrary $s$ the $s=2$ case considered in \cite{GG1}, we define:
\begin{dfn}
\label{dfn:clif}
For $p\in X$, the Clifford algebra $\mathcal{C}_p$ is 
\[\mathcal{C}_p = Cl(T_pX \oplus T^*_pX,q)\]
with the quadratic form $q$ induced by the  metric
\[\begin{array}{ll}\forall i,j,h,k & <v_{ij},v_{hk}> = 0\\
\forall i,j,h,k & <\frac{\partial}{\partial v_{ij}},\frac{\partial}{\partial v_{hk}}> = 0\\
\forall (i,j)\not= (h,k) &  <v_{ij},\frac{\partial}{\partial v_{hk}}> = 0\\
\forall i,j & <v_{ij},\frac{\partial}{\partial v_{ij}}> = -\frac{1}{2}
\end{array}\]
\end{dfn}
\begin{rmk} The Clifford algebras $\mathcal{C}_p$ for varying $p$ define a Clifford bundle $\mathcal{C}$ on $X$, as the definition of $\mathcal{C}_p$ is independent on the choice of a basis. Indeed, the quadratic form used to define it is simply induced by $-\frac{1}{2}$ times the natural bilinear pairing  $T_pX \otimes T^*_pX\to \R$.
\end{rmk}
\begin{pro}
The Clifford algebra $\mathcal{C}_p$ has a canonical representation $\rho_p$ on $\bigwedge T^*_p X$, induced by the wedge and contraction operators $E_{ij}$ and $I_{ij}$ via the map
\[\rho_p(v_{ij}) =  E_{ij},\qquad \rho_p\left(\frac{\partial}{\partial v_{ij}}\right) =  I_{ij}\]
\end{pro}
\dimo The Clifford relations
\[\phi\psi + \psi\phi =  -2<\phi,\psi>\]
are precisely the content of Proposition \ref{pro:clifrelations}.
\qed\\
Abusing slightly the notation, we will identify $\mathcal{C}_p$ with its (faithful) image inside $End_\R\left(\bigwedge^*T^*_pX\right)$, and we will omit any reference to the map $\rho_p$ when it will not be necessary. Actually, as the representation above is a real analogue of the Spinor representation, it is easy to check that the map $\rho_p$  is an isomorphism of associative algebras. One then has:
\begin{dfn}
\label{dfn:quadratic} The linear subspace $\mathcal{C}^2_p$ of $\mathcal{C}_p$ is the image of the natural map $\bigwedge^2(T_pX\oplus T_p^*X)\to \mathcal{C}_p$. The linear subspace $\mathcal{C}^0_p$ of $\mathcal{C}_p$ is the subspace generated by $1$. 
\end{dfn}
Recall (see for instance \cite{LM}) that $\mathcal{C}_p^2$ is a Lie subalgebra  of $\mathcal{C}_p$ (with the commutator bracket).
\begin{pro} The bundle of Lie algebras $\ggr$ is a sub-bundle of  $\mathcal{C}^2$. Any local determination of the operator $J$ is a (local) section of $\mathcal{C}^2$.
\label{pro:subclif}
\end{pro}
\dimo Let us fix $p\in X$. We consider the pointwise values of the operators $L_{\alpha\beta}$, the $V_j$ and the $A_j$; they all lie inside $\mathcal{C}_p^2\oplus \mathcal{C}_p^0$ by Proposition \ref{pro:clifrelations} and by the fact that the forms $\omega_{ij}$ restrict to elements of $\bigwedge^2T^*_pX$.  The space $<J>$ lies inside $\mathcal{C}_p^2\oplus \mathcal{C}_p^0$ by Proposition \ref{pro:dfnJ}. By definition the elements $\mathcal{C}_p^2$  are commutators, and therefore have trace zero in any representation, and hence also in the $\rho_p$. Moreover, again by inspection all the generators of the fibre in $p$ of $\ggr$ have trace zero once represented via $\rho_p$ (they are nilpotent), and therefore they must lie inside $\mathcal{C}_p^2$. Both pointwise determinations of operator $J$ are in the Lie algebra of the isometry group, and therefore they too have trace zero and hence sit inside $\mathcal{C}_p^2$. As $\mathcal{C}_p^2$ is closed under the commutator bracket of $\mathcal{C}_p$, and this commutator coincides with the composition bracket of operators, we have the conclusion.
\qed
\begin{rmk} For any $p\in X$, the Clifford algebra $\mathcal{C}_p$ is isomorphic to the standard Clifford Algebra $\mathbf{Cl}_{2s+2,2s+2}$, as the metric used to define it has signature $(2s+2,2s+2)$. The previous proposition therefore shows that all the fibres of $\ggr$ are Lie subalgebras of $\mathbf{Cl}_{2s+2,2s+2}^2 \cong \mathbf{spin}_{2s+2,2s+2}$.
\end{rmk}
\begin{rmk}
For any fixed $p\in X$, giving degree $1$ to the operators $E_{ij}$ and degree $-1$ to the operators $I_{ij}$, we induce a $\Z$-degree on  $\mathcal{C}_p$. This degree coincides with the degree of the operators induced from the grading on  the forms from $\bigwedge^*T^*_pX$.
\end{rmk}
Similarly to Definition \ref{dfn:clif}, for any $p\in X$ one could define a Clifford Algebra 
\[\mathcal{C}l(\R^{s+1}\oplus (\R^{s+1})^*,q_{nat})\]
 where $q_{nat}$ is the quadratic form induced by ($-\frac{1}{2}$ times) the natural paring and 
\[\R^{s+1} = <\tilde{E}_0,...,\tilde{E}_s>,\quad  (\R^{s+1})^* = <\tilde{I}_0,...,\tilde{I}_s>\]
One has also a natural representation on $\bigwedge^*T_pX$ of the operators $[{\tilde E}_j,\tilde{E}_k]$, 
$[\tilde{E}_j,\tilde{I}_k]$,
$[\tilde{I}_j,\tilde{I}_k]$ 
generating the  degree two part of this Clifford Algebra , induced by the map which acts as follows:
\[
\begin{array}{lcl}
{[\tilde{E}_j,\tilde{E}_k]} & \to & 2L_{jk}\\
{[\tilde{E}_j,\tilde{I}_k]} & \to & 2L_{j\bk}\\
{[\tilde{I}_j,\tilde{I}_k]} & \to & 2L_{\bj\bk}
\end{array}
\]
This gives directly the bundle $\lefs$  as a quotient of the $\mathbf{spin}_{s+1,s+1}$ Lie Algebra bundle of this Clifford bundle, proving again that its fibre is indeed $\mathbf{so}(s+1,s+1,\R)$.

\section{Quadratic invariants and Hodge decomposition}
\label{sec:plethysm}
Fixing $p\in X$ and a determination $J$ at $p$, the complex structure $J$ acts on all the Clifford algebra $\mathcal{C}_p$ by adjunction with respect to the commutator bracket, and sends its quadratic part $\mathcal{C}^2_p$  to itself from Proposition \ref{pro:subclif}. 
\begin{dfn} We call {\em quadratic invariants} the elements in $\mathcal{C}^2_p$ which commute with $J$. For varying $p$, we obtain a bundle of quadratic invariants.
\end{dfn}

As usual, to decompose the representation $\bigwedge^* T^*X$ with respect to the weight induced by $J$, it is necessary to consider complexified forms (and algebras). The weight decomposition of the  space $T^*_p X \otimes \C$ is obtained introducing a new basis for each $W_j\otimes\C =
<v_{1j},v_{2j}>_\C$:
\[w_j = \frac{1}{\sqrt{2}}(v_{1j} + \imath\ v_{2j}),\qquad\overline{w}_j = \frac{1}{\sqrt{2}}( v_{1j} - \imath
v_{2j})\]
To describe  explicitely the space of (complex) quadratic invariants in the Clifford algebra $\mathcal{C}_p$, let us introduce the following notation, which gives a basis of eigenvectors for the (adjoint) action of $J$:
\begin{dfn}
\[E_{w_j} =  \frac{1}{\sqrt{2}}(E_{1j} + \imath E_{2j}),\qquad E_{\overline{w}_j} =  \frac{1}{\sqrt{2}}(E_{1j} - \imath E_{2j})\]
\[I_{w_j} = \frac{1}{\sqrt{2}}( I_{1j} - \imath I_{2j}),\qquad I_{\overline{w}_j} = \frac{1}{\sqrt{2}}( I_{1j} + \imath I_{2j})\] 
\end{dfn}

\begin{lem} 
\label{lem:weights} The adjoint action of the complex structure operator  $J$  on $E_{w_j}$, $I_{w_j}$, $E_{\overline{w}_j}$, $I_{\overline{w}_j}$ is:
\[[J,E_{w_j}] = -\imath E_{w_j},\quad [J,I_{w_j}] = \imath I_{w_j}\]
\[[J,E_{\overline{w}_j}] = \imath E_{\overline{w}_j},\quad[J,I_{\overline{w}_j}] = -\imath I_{\overline{w}_j}\]
\end{lem}
\dimo It is enough to consider the corresponding $J$-weights of the $w_j,\overline{w}_j$. \qed\\
As $\gg\subset \mathcal{C}^2\otimes\C$ from Proposition \ref{pro:subclif}, in the following we show that $\gg$ lies inside the bundle of quadratic invariants. Immediately after we will give a basis for the space of quadratic invariants,  thus providing a first  upper bound for $\gg$ (which will be later shown to be off by only 1).
\begin{pro}
\label{pro:ggrinvariant}
The operator $J$ commutes with all the elements in the fiber at $p$ of $\gg$. 
\end{pro}
\dimo We prove the statements by a direct computation. It is useful to rewrite $\omega_{jk}$
(and hence $L_{jk}$ which is wedge with $\omega_{jk}$) in terms of the
basis generated by the $w_j,\overline{w}_k$:
\[\omega_{jk} = v_{1j}\wedge v_{1k} + v_{2j}\wedge v_{2k} =
\frac{1}{2}\left(w_{j}\wedge \overline{w}_{k} + 
\overline{w}_{j}\wedge w_k\right)\] and therefore 
\[L_{jk} = \frac{1}{2}\left([E_{w_j},E_{\overline{w}_k}] - [E_{w_k},E_{\overline{w}_j}]\right)\]
The vanishing $[J,L_{jk}] = 0$ then follows immediately from Lemma \ref{lem:weights}.
Similarly, to show that $[J,V_k] = 0$ it is enough to observe that
\[V_k =  \frac{\imath}{2}[ E_{w_{k}}, E_{\overline{w}_{k}}]\] 
The corresponding commutation relations for the adjoint
operators follow from the fact that $J^*
= -J$, as noticed in Remark \ref{rmk:dfnJ}. \qed \\

The following proposition will show that, except for a toral part which will be discussed later, all the quadratic invariants of the Clifford bundle $\mathcal{C}$ lie inside $\ggrs\subset \gg$. It will follow therefore that
\[4 (s+1)^2-2(s+1) \leq dim_\R\, \ggrs \leq dim_\C \,\gg \leq 4 (s+1)^2\]
\begin{pro}
The following  \(4 (s+1)^2\) operators are  a linear basis for the quadratic $J$-invariants: 
\begin{enumerate} 
\item $[E_{w_i},E_{\overline{w}_j}]$ with \(i\neq j\).
\item $[I_{w_i},I_{\overline{w}_j}]$ with \(i\neq j\).
\item $[E_{w_i},E_{\overline{w}_i}]$ where \(i=0,1,\ldots , s\). 
\item $[I_{w_i},I_{\overline{w}_i}]$ where \(i=0,1,\ldots , s\). 
\item $[E_{w_i},I_{w_j}]$ with \(i\neq j\).
\item $[E_{\overline{w}_i},I_{\overline{w}_j}]$ with \(i\neq j\).
\item $[E_{w_i},I_{w_i}]$ where \(i=0,1,\ldots , s\). 
\item $[E_{\overline{w}_i},I_{\overline{w}_i}]$ where \(i=0,1,\ldots , s\).
\end{enumerate}
The  \(4 (s+1)^2-2(s+1)\) operators of type  \((1),(2),(3),(4),(5),(6)\) belong to the bundle of real algebras  \(\ggrs\subset \gg\).
\label{pro:basis}
\end{pro}
\dimo In this proof, we fix $p\in X$ and all the bundles and operators will be considered at this point. The $J$-weight of a bracket of $J$-homogeneous operators is the sum of the respective weights. The quadratic "monomials" (with respect to the bracket) in the $E_{w_j},I_{w_j},E_{\overline{w}_j},I_{\overline{w}_j}$ are all $J$-homogeneous, and therefore to find a basis of $J$-invariant quadratic operators it is enough to identify the $J$-invariant quadratic monomials. To be $J$-invariant means simply to have weight zero, and the computation of the $J$-weight of the quadratic mononials follows immediately from those of $E_{w_j},I_{w_j},E_{\overline{w}_j},I_{\overline{w}_j}$, which are respectively $-\imath,\imath,\imath,-\imath$. 

It remains to be shown that the  monomials of type \((1),(2),(3),(4),(5),(6)\)  belong to \(\ggrs\).
Let us consider 
$[E_{w_i},E_{\overline{w}_j}]$ with \(i\neq j\). Since $E_{w_i}$ and $E_{\overline{w}_j}$ anticommute, this is equal to $2 E_{w_i}E_{\overline{w}_j}$.
Then
\[2E_{w_i} E_{\overline{w}_j}=(E_{1i}+\imath E_{2i})(E_{1j}-\imath E_{2j})=E_{1i}E_{1j}+E_{2i}E_{2j}+\imath (E_{2i}E_{1j}-E_{1i}E_{2j})=\] \[=L_{ij}+\imath (E_{2i}E_{1j}-E_{1i}E_{2j})\]

We have therefore to show that \(\imath (E_{2i}E_{1j}-E_{1i}E_{2j})\) belongs to \(\ggrs\).

We recall that, by Corollary \ref{cor:ibj}, the elements \(L_{i\bj}\) belong to \(\ggrs\) and notice that
\[[L_{i\bj}, \imath V_j]=\imath (E_{1i}I_{1j}+E_{2i}I_{2j})E_{1j}E_{2j}-\imath E_{1j}E_{2j}(E_{1i}I_{1j}+E_{2i}I_{2j})= \] \[=\imath (E_{1i}E_{2j}-E_{1i}E_{1j}I_{1j}E_{2j}-E_{2i}E_{1j}+E_{2i}E_{1j}E_{2j}I_{2j}-E_{1j}E_{2j}E_{1i}I_{1j}-E_{1j}E_{2j}E_{2i}I_{2j})=\]\[= \imath (E_{1i}E_{2j}-E_{2i}E_{1j})\]
which concludes the proof that the monomial $[E_{w_i},E_{\overline{w}_j}]$ lies in \(\ggrs\).
By adjunction, we immediately have that also $[I_{w_i},I_{\overline{w}_j}]$ lies in \(\ggrs\).

Also the monomials $[E_{w_i},E_{\overline{w}_i}]$ belong to \(\ggrs\); in fact they are immaginary multiples of the volume forms: 
\[2[E_{w_i}E_{\overline{w}_i}]=(E_{1i}+\imath E_{2i})(E_{1i}-\imath E_{2i})-(E_{1i}-\imath E_{2i})(E_{1i}+\imath E_{2i})=\]\[= 2 \imath (E_{2i}E_{1i}-E_{1i}E_{2i})=-4\imath V_i\] 
As a consequence, we have, by adjunction, $2[I_{w_i},I_{\overline{w}_i}]=-4\imath A_i$.

Let us consider the monomials of type (5)
$[E_{w_i},I_{w_j}]$ with \(i\neq j\). Since $E_{w_i}$ and $I_{w_j}$ anticommute, this  is equal to $2 E_{w_i}I_{w_j}$.
Then
\[2E_{w_i}I_{w_j}=(E_{1i}+\imath E_{2i})(I_{1j}-\imath I_{2j})=E_{1i}I_{1j}+E_{2i}I_{2j}+\imath (E_{2i}I_{1j}-E_{1i}I_{2j})=\] \[=L_{i\bj}+\imath (E_{2i}I_{1j}-E_{1i}I_{2j})\]

We have therefore to show that \(\imath (E_{2i}I_{1j}-E_{1i}I_{2j})\) belongs to \(\ggrs\).

We notice that
\[[L_{\bi\bj}, \imath V_j]=\imath (I_{1i}I_{1j}+I_{2i}I_{2j})E_{1j}E_{2j}-\imath E_{1j}E_{2j}(I_{1i}I_{1j}+I_{2i}I_{2j})= \] \[=\imath (I_{1i}E_{2j}-I_{1i}E_{1j}I_{1j}E_{2j}-I_{2i}E_{1j}+I_{2i}E_{1j}E_{2j}I_{2j}-E_{1j}E_{2j}I_{1i}I_{1j}-E_{1j}E_{2j}I_{2i}I_{2j})=\]\[= \imath (I_{1i}E_{2j}-I_{2i}E_{1j})\]
which by adjunction gives that also 
\[[L_{ji},\imath  A_j]=\imath(E_{1i}I_{2j}-E_{2i}I_{1j})\]
This allows us to conclude that the monomial $[E_{w_i},I_{w_j}]$ (as of course its conjugate $[E_{\overline{w}_i},I_{\overline{w}_j}]$) lies in \(\ggrs\).\qed \\




\section{All the fibres of the bundle \(\gg\) are isomorphic to \(\mathbf{sl}(2s+2,\C)\).}
\label{sec:isoc}
In this section, we fix once and for all a determination of $J$ at the point $p$ and consider the Hodge decomposition of $\bigwedge_\C ^*T_p^* X$ with respect to the (almost) complex structure $J$. We will use this information to first study the complex algebra $\gg$, while in the next sections we will concentrate on its reals forms. In the rank $3$ case this corresponds to performing the plethysm with respect to the action of $\mathbf{SO}(3,\R)$ (see \cite{GG2}, where we analyze the case $s=2$ in the context of WSD manifolds).
The Hodge (type) decomposition of forms on $X$ with respect to the complex structure $J$ 
\[{\bigwedge^k}_\C T_p^* X = \bigoplus_{r+t = k}{\bigwedge^{r,t}}_\C T^*_p X\]
can be described as usual explicitely as follows, using the $J$-homogeneous basis $w_j,\overline{w}_k$:
\[\bigwedge^{r,t}T^*_\C X_p = <w_{i_1}\wedge\cdots\wedge w_{i_r}\wedge
\overline{w}_{j_1}\wedge\cdots\wedge
\overline{w}_{j_t}~|~i_1,...,j_t~\in\{0,1,2,\ldots ,s\}>_\C\] 
\begin{dfn}
\label{dfn:ialpha} 
At a given point $p\in X$, and with the chosen a determination of $J$ at $p$,
we indicate with $\mathcal{I}_\alpha$ the subspace (isotypical component)  of forms of $J$-weight $\imath\alpha$ ($-s-1\leq \alpha \leq s+1$).
\end{dfn}
Here is for instance the Hodge "diamond" in the case $s=4$ (we used the following notation:  the symbol  \(\Omega^{r,t}_m\) indicates the space \(\bigwedge^{r,t}T^*_\C X_p = \Omega^{r,t}_\C X_p\)  and specifies its dimension \(m\)). The Lie algebra bundle $\gg$ acts preserving the weight of forms, and therefore the spaces $\mathcal{I}_\alpha$ which are the columns in the Hodge diamond.
{\tiny
\begin{table}[h]
\begin{center}
\begin{tabular}{|c|c|c|c|c|c|c|c|c|c|c|c|}
\hline
 & &  &       &  & &  &       & & &  & \\
   & \({\mathcal I}_{-5}\)  & \({\mathcal I}_{-4}\)  & \({\mathcal I}_{-3}\)  &\( \, \)\({\mathcal I}_{-2}\)  \( \, \) & \( \, \)\({\mathcal I}_{-1}\)  \( \, \)  &\( \, \) \({\mathcal I}_{0}\) \( \, \)&\( \, \)\({\mathcal I}_{1}\) \( \, \)&\( \, \)\({\mathcal I}_{2}\)  \( \, \) & \({\mathcal I}_{3}\)  &\({\mathcal I}_{4}\)   & \({\mathcal I}_{5}\)    \\
 & &  &       &  & &  &       & & &  & \\
\hline
 & &  &       &  & &  &       & & &  & \\
 $\Omega^0_\C X_p$  & &  &       &  & &\(\Omega^{0,0}_1\)  &       & & &  & \\
 & &  &       &  & &  &       & & &  & \\
 & &  &       &  & &  &       & & &  & \\
$\Omega^1_\C X_p$    & &  &       &  &\(\Omega^{1,0}_5\)  &  & \(\Omega^{0,1}_{5}\)      & & &  & \\
 & &  &       &  & &  &       & & &  & \\
 & &  &       &  & &  &       & & &  & \\
  $\Omega^2_\C X_p$    & &  &       & \(\Omega^{2,0}_{\binom{5}{2}} \)& & \(\Omega^{1,1}_{5^2}\) &       &\(\Omega^{0,2}_{ \binom{5}{2}}\) & &  & \\
 & &  &       &  & &  &       & & &  & \\
 & &  &       &  & &  &       & & &  & \\
 $\Omega^3_\C X_p$    & &  &   \(\Omega^{3,0}_{\binom{5}{3}} \)    &   & \(\Omega^{2,1}_{5\binom{5}{2}} \)  &    & \(\Omega^{1,2}_{5\binom{5}{2}} \)         &    &   \(\Omega^{0,3}_{\binom{5}{3}} \)  &  & \\
 & &  &       &  & &  &       & & &  & \\
 & &  &       &  & &  &       & & &  & \\
  $\Omega^4_\C X_p$    & &  \(\Omega^{4,0}_{5}\)  &       & \(\Omega^{3,1}_{\binom{5}{4}\binom{5}{2}} \)& & \(\Omega^{2,2}_{\binom{5}{2}\binom{5}{2}} \) &  & \(\Omega^{1,3}_{\binom{5}{4}\binom{5}{2}} \)      &   &  \(\Omega^{0,4}_{5}\)  & \\
 & &  &       &  & &  &       & & &  & \\
 & &  &       &  & &  &       & & &  & \\
  $\Omega^{5}_\C X_p$    & \(\Omega^{5,0}_{1}\)  &  &  \(\Omega^{4,1}_{5^2}\)    &  & \(\Omega^{3,2}_{\binom{5}{3}\binom{5}{2}}\) &  & \(\Omega^{2,3}_{\binom{5}{3}\binom{5}{2}}\)    & & \(\Omega^{1,4}_{5^2}\) &  & \(\Omega^{0,5}_{1}\)  \\
 & &  &       &  & &  &       & & &  & \\
  & &  &       &  & &  &       & & &  & \\
  $\Omega^{6}_\C X_p$     & &  \(\Omega^{5,1}_{5}\)  &       & \(\Omega^{4,2}_{\binom{5}{4}\binom{5}{2}} \)& & \(\Omega^{3,3}_{\binom{5}{2}\binom{5}{2}} \) &  & \(\Omega^{2,4}_{\binom{5}{4}\binom{5}{2}} \)      &   &  \(\Omega^{1,5}_{5}\) & \\
 & &  &       &  & &  &       & & &  & \\
  & &  &       &  & &  &       & & &  & \\
 $\Omega^{7}_\C X_p$      & &  &   \(\Omega^{5,2}_{\binom{5}{3}} \)    &   & \(\Omega^{4,3}_{5\binom{5}{2}} \)  &    & \(\Omega^{3,4}_{5\binom{5}{2}} \)         &    &   \(\Omega^{2,5}_{\binom{5}{3}} \)  &  & \\
  & &  &       &  & &  &       & & &  & \\
     $\Omega^{8}_\C X_p$      & &  &       & \(\Omega^{5,3}_{ \binom{5}{2}} \)& & \(\Omega^{4,4}_{5^2}\) &       &\(\Omega^{3,5}_{ \binom{5}{2}} \) & &  & \\
      & &  &       &  & &  &       & & &  & \\
  & &  &       &  & &  &       & & &  & \\
     $\Omega^{9}_\C X_p$     & &  &       &  &\(\Omega^{5,4}_{5}\) &  & \(\Omega^{4,5}_{5}  \)    & & &  & \\
 & &  &       &  & &  &       & & &  & \\
  & &  &       &  & &  &       & & &  & \\
   $\Omega^{10}_\C X_p$  & &  &       &  & &\(\Omega^{5,5}_{1}\)  &       & & &  & \\
 & &  &       &  & &  &       & & &  & \\
 & &  &       &  & &  &       & & &  & \\
\hline
 & &  &       &  & &  &       & & &  & \\
  Dim.   & 1  & \(10\)  &\( \binom{10}{2}   \) &\( \binom{10}{3}   \) &  \( \binom{10}{4}  \)&   \(\binom{10}{5}   \)    & \( \binom{10}{4}  \) & \( \binom{10}{3}     \) &\( \binom{10}{2}   \) &\(10\)   & 1 \\
 & &  &       &  & &  &       & & &  & \\
\hline
\end{tabular}
\medskip
\caption{Hodge diamond (case \(s=4\))} \label{tab:isotyp}
\end{center}
\end{table}

}

\begin{teo}
\label{teo:slstructure} 
Let $X$ be a (almost, pointwise) $s$-\ka manifold or rank two. \medskip\\
a) ~
The Lie algebra bundle \(\gg\) has fibre isomorphic to \(\mathbf{sl}(2s+2,\C)\).\\
b)~ At a given point $p\in X$, the direct sum of ${\gg}_{,p}$ with the space spanned by the operator $J_p$ is the set of all quadratic invariants of $\mathcal{C}_p$.\\
c)~At a given point $p\in X$, the restriction of $\gg$ to the $2s +2$ dimensional space $\mathcal{I}_{-s}$ of forms of $J$-weight $-s$ is faithful.
\end{teo}

\dimo We work at a fixed point $p$.
The isotypical    component   \({\mathcal I}_{-s}\)     has dimension \(2s +2\) and has a basis \(\{b_i\}\) (\(0\leq i\leq 2s +2\)) given by the following monomials:\\

\begin{itemize}
\item  \(b_i=w_{0}\wedge \ldots \wedge \widehat {w_i} \wedge \cdots\wedge w_{s}\), for \(i \in \{0,1,2,\ldots , s\}\),  where \( \widehat {w_i} \) means that  \(w_i\) is omitted and therefore the monomial  has degree \(s\).
\item  \(b_{s+1+i}=w_{0}\wedge \cdots\wedge w_{s}\wedge \overline{w}_i\), where \(i \in \{0,1,2,\ldots , s\}\) and the monomial has degree \(s+2\).

\end{itemize}
\bigskip

It is then immediate to check that, for instance, 
 \[[E_{w_0},I_{w_1}](b_0)= [E_{w_0},I_{w_1}](w_{1}\wedge \cdots\wedge w_{s})=2w_{0}\wedge \widehat {w_1} \wedge \cdots\wedge w_{s}=2b_1\qquad \; \;  \; \;   \; \;  \]
 \[[E_{w_s},E_{\overline{w}_0}](b_s)=[E_{w_s},E_{\overline{w}_0}](w_{0}\wedge \cdots\wedge w_{s-1})=2w_{0}\cdots \wedge w_{s} \wedge \overline{w}_0=2b_{s+1} \; \; \;  \; \;   \, \]
 \[[E_{\overline{w}_1},I_{\overline{w}_0}](b_{s+1})=[E_{\overline{w}_1},I_{\overline{w}_0}](w_{0}\wedge \cdots \wedge w_{s} \wedge \overline{w}_0 )=2w_{0}\cdots \wedge w_{s} \wedge \overline{w}_1=2b_{s+2} \]

Completely analogous computations  show that, when we represent the action of \(\ggp\) on  the isotypical component  \({\mathcal I}_{-s}\)  using the above mentioned basis,
all the elementary matrices \(e_{ij}\) (where \(i\neq j\) and \(e_{ij} \) is the matrix with all the entries equal to 0 except for the entry \((i,j)\) which is  1) are obtained using the quadratic invariants of type    \((1),(2),(3),(4),(5),(6)\) which in Proposition \ref{pro:basis} were shown to lie in \(\ggp\).

More precisely, we have the following identifications for the ``positive'' set of Serre generators \(e_{j+1,j}\):\\
\begin{itemize}
\item  \(e_{j+2,j+1}=\frac{1}{2}[E_{w_j},I_{w_{j+1}}]\) for \(0\leq j\leq s-1\); 
\item \(e_{s+2,s+1}=\frac{1}{2}[E_{w_s},E_{\overline{w}_0}]\);
\item  \(e_{s+3+j,s+2+j}=\frac{1}{2}[E_{\overline{w}_{j+1}},I_{\overline{w}_j}]\) for \(0\leq j\leq s-1\).
\end{itemize}
\bigskip

Therefore \(\ggp \) acts as \(\mathbf{sl}(2s+2,\C)\) on   \({\mathcal I}_{-s}\)   (notice that, as  the generators \(L_{ij}\), \(\Lambda_{ij}=L_{\bj\bi}\), \(V_i\), \(A_i\) of \(\ggp\) are nilpotent, they still have trace zero when restricted to     \({\mathcal I}_{-s}\) ).

Summing up, the algebra \(\ggp\) has a quotient isomorphic to the simple algebra \(\mathbf{sl}(2s+2,\C)\) and is embedded in the  \(4 (s+1)^2\)-dimensional space of the quadratic invariants; now, since  the quadratic invariant \(J_p\) doesn't belong to \(\ggp\) (in fact the trace of  its restriction  to \({\mathcal I}_{-s}\)  is different from 0, since \(J_p\) acts on  \({\mathcal I}_{-s}\)  as multiplication by \(-\imath s\)), we conclude that \(\gg\) has dimension \(4 (s+1)^2-1\). Therefore the restriction to       \({\mathcal I}_{-s}\)  provides us with an  isomorphism of \(\ggp\) with   \(\mathbf{sl}(2s+2,\C)\).
\qed\\

The decomposition $T^*X = W_0 \oplus W_1\oplus\cdots \oplus
W_s$ induces naturally a multi-degree on $\bigwedge^*T^*_\C X$ with
values in $\Z^{s+1}$, which we indicate with $mdeg$. This follows from
the equation
\[\bigwedge^n T_\C^*X \cong \bigoplus_{p_0+p_1+\cdots +p_s = n}
\bigwedge^{p_0}\left(W_0\otimes\C\right)\oplus
\bigwedge^{p_1}\left(W_1\otimes\C\right)\oplus\cdots \oplus 
\bigwedge^{p_s}\left(W_s\otimes\C\right)\] We notice furthermore that
the (complexified) decomposition above is preserved by the operator
$J$, and therefore $mdeg$ commutes with the action of $\so2$. 
We will still call by \(mdeg\) the multidegree induced on  the bundle \(\gg\) by the previous one.

\begin{cor}
For any fixed point $p\in X$, the space of quadratic invariants with \(mdeg\) equal to \((0,0,0,\ldots ,0)\) coincides with the \(2s +2\)-dimensional space spanned by  the quadratic monomials  $[E_{w_i},I_{w_i}]$ and by  their   complex conjugates $[E_{\overline{w}_i},I_{\overline{w}_i}]$. This space  can be expressed as   the direct sum of  a maximal toral subalgebra of  \({\gg}_{,p}\cong  \mathbf{sl}(2s +2,\C)\) plus the one dimensional  subspace \(<J_p>\).

\end{cor}

\dimo The first assertion is trivial (by inspection of the basis of the quadratic invariants described in Proposition \ref{pro:basis}). The restriction   of \(\ggp\) to  \({\mathcal I}_{-s}\) is an isomorphism by the previous theorem. Therefore, the diagonal matrices with trace zero (in the same basis $\{b_i\}$ used in the proof of the theorem) provide a toral subalgebra of $\ggp$ formed by operators with $mdeg$ equal to zero as they can be obtained as brackets of operators with opposite $mdeg$.
As the above mentioned basis is made up of \(mdeg\)-homogeneous elements, all the quadratic invariants with vanishing \(mdeg\) must be   associated  to diagonal matrices. Summing up, the quadratic invariants with vanishing $mdeg$ are  the toral elements of $\ggp$ plus $J_p$ (which has certainly vanishing $mdeg$, since it admits a global basis of $mdeg$-homogeneous eigenvectors).\qed \\

\section{The real forms of $\gg$}
In the previous section we described completely the complex bundle of Lie algebras $\gg$ (see Theorem \ref{teo:slstructure}). Here, we will give a complete description or the geometrically natural bundles of real Lie algebras $\ggrs$ $\ggru$ and $\ggr$.\\
Recall that the bundle of Lie algebras $\ggrs$ is generated by the Lefschetz operators $L_{\alpha\beta}$ and by the $\imath V_j,\imath A_k$.
\label{sec:split}
\begin{teo}
\label{pro:realform} The bundle of real Lie algebras \(\ggrs\) has fibre isomorphic to the Lie algebra \(\mathbf{sl}(2s +2,\R)\).

\end{teo}

\dimo  We identify the fibre of \(\gg\) at a point $p$ with \(\mathbf{sl}(2s +2,\C) \) using the faithful representation \({\mathcal I}_{-s}\) and the basis $\{b_i\}$ of  \({\mathcal I}_{-s}\)  provided in the proof of Theorem  \ref{teo:slstructure}.  As we already noticed in that proof, in this basis all the invariant monomials of type  \((1),(2),(3),(4),(5),(6)\) act via real matrices, and provide (up to a scalar) all the elementary matrices \(e_{ij}\) (\(i\neq j\)). Therefore they generate over the real numbers the  subalgebra \(\mathbf{sl}(2s +2,\R) \) of \(\mathbf{sl}(2s +2,\C) \). 
If we prove  that the whole fibre at $p$ of \(\ggrs\) lies in \(\mathbf{sl}(2s +2,\R) \), then,  since \(\ggrs\) is a real form of \(\gg\), we must have \({\ggrs}_{p}= \mathbf{sl}(2s +2,\R) \).

It suffices to notice that the generators \(L_{ij}\) and \(\imath V_j\) (and therefore  their adjoints) are in the real algebra generated by the above mentioned  invariant monomials.
Now, in the proof of Theorem  \ref{teo:slstructure},  \(\imath V_j\) was obtained  as \(\displaystyle{-\frac{1}{2} [E_{w_i},E_{\overline{w}_i}]}\) and in Proposition \ref{pro:ggrinvariant} we showed that 
\[ [E_{w_i},E_{\overline{w}_j}]-[E_{w_j},E_{\overline{w}_i}]=2L_{ij}\]
\qed\\

On the complex bundle of vector spaces $\bigwedge^*_\C T^*X$ there is a natural hermitean inner product $<~,~>$, obtained from the wedge operation on forms (cf. \cite{GG2} where we used a superHermitean variant of this product for the rank $3$ case), and defined below. Associated to this pairing, there is a natural notion of antihermitean operator. We will prove that the set of antihermitean operators inside $\gg$ is a real form for $\gg$, generated by operators naturally derived from the geometry and coinciding with $\ggru$. 

\begin{dfn}
\label{dfn:hermitean}
For every $p\in X$ there is a natural  non degenerate Hermitean inner product $<~,~>_p$ on $\bigwedge_\C^*T^*_pX$, defined using the natural (standard) Hermitean inner product $(~,~)_p$ associated to the metric $\me$ and the (pointwise) volume form $\Omega$ associated to the metric $\me$ and to the chosen determination of $J$ at $p$:
\[ <\alpha,\beta>_p = \imath^{deg(\alpha) deg(\beta)}(\alpha\wedge \overline{\beta},\Omega)_p\]
We indicate with $<~,~>$ the corresponding form with values in smooth functions.
\end{dfn}
Let us denote with $\ns$ the (complex linear) operator obtained composing conjugation with the Hodge star associated to the metric. 
\begin{pro} For every $p\in X$, the pairing $<~,~>_p$ satisfies the following properties:\\
\label{primeprop}

a) $\quad < \alpha,\beta>_p = \imath^{(deg({\alpha}) +2)deg(\beta))}(\alpha, \ns\beta)_p$\\

b) $\quad <~,~>_p$ is  preserved by the the operator $J$ in derived sense, namely 
\[\forall\alpha\beta~ <J\alpha,\beta> + <\alpha,J\beta> = 0\]

c) $\quad <~,~>_p$ is  preserved by the operator $\ns$, namely 
\[\forall\alpha\beta~ <\ns\alpha,\ns\beta> = <\alpha,\beta> \]

d) The pure weight components $\mathcal{I}_k$ are mutually $\quad <~,~>_p$-orthogonal and  $\quad <~,~>_p$ is nondegenerate when restricted to any one of them.
\end{pro}

\dimo The first three facts are standard. For the orthogonality statement in part $d)$, we observe that,  if $\alpha\in\mathcal{I}_h$ and $\beta\in\mathcal{I}_k$ with $deg(\alpha)+ deg(\beta) = dim(X)$ then 
\[<\alpha,\beta>_p\Omega = (\alpha\wedge\overline{\beta})_p\]
is a complex number times a form of $J$-weight zero, but from the right hand side it also must have $J$-weight equal to $(h-k)\imath$. Therefore if $h\not= k$, it must be zero.\\
Restricting to a single $\mathcal{I}_k$, notice that $\ns$ sends this component to itself (as it commutes with $J$), and then if $\alpha\not=0$ in $\mathcal{I}_k$, 
$<\alpha,\ns\alpha> $ is a power of $\imath$ times  $(\alpha,\alpha)$ by point $a)$, and is therefore nonzero.
\qed\\

We want now to characterize the operators inside $\gg_{,p}$ which preserve the form $<~,~>_p$. 
First we observe that, since the dimension of \(T^*_pX\) is even,  \(\ns\ns\) is equal to the identity on the forms of even degree  while \(\ns\ns=-I\) when restricted to the odd forms. Then, for fixed $p\in X$, using the expression $<\alpha,\beta>_p  = (\alpha,\ns\beta)_p$,   we see that the ``differential'' condition for preservation of the form by the operator $\phi$
\[\forall\alpha\forall\beta \quad <\phi(\alpha),\beta>_p + <\alpha,\phi(\beta)>_p = 0\]
is equivalent to $\phi^\ss = -\ns\,\phi\,\ns$ on the even forms and to  $\phi^\ss = \ns\,\phi\,\ns$ on the odd forms.\\

The next two theorems show that the bundle $\ggru$ (generated at any point by the value of  the operators    \(\imath L_{j,k} \) (\(j\neq k\)), \(\imath V_i\) and  their adjoints, see Definition \ref{dfn:ggru}) is  precisely the bundle of Lie subalgebras given point by point by the operators which preserve the form $<~,~>$:

\begin{teo}
\label{teo:antihermitean} The Lie algebra bundle \(\ggru\)  preserves the  form $<~,~>$. 
\end{teo}
\dimo
As we observed before, the statement is equivalent to the fact that the condition $\phi^\ss (\alpha) = (-1)^{deg \ \alpha +1} \ns\,\phi\,\ns (\alpha) $ holds for all the sections $\phi$ of $\ggru$ and all the homogeneous elements \(\alpha \in \bigwedge_\C^*T^*X\).\\
It is enough to check these equations for the generators of $\ggru$, at a fixed point $p\in X$, which as usual we omit from the notation for the operators when not strictly necessary.

Let \(\Psi\) be one of  the generators $\imath L_{ij}$ or \(\imath V_k\);  this means that  \(\Psi\) is the operator given by the wedge with an even form \(\imath \psi\), where \(\psi\) is real. One has, given \(v\) homogeneous of degree $h$ and \(w\) in \(\bigwedge_\C^*T^*_pX\) with degree of the same parity (which is the only possibly non-vanishing case):
\[(\Psi (v), w)_p=
(\imath \psi \wedge v \wedge \overline{ \ns w},\Omega)_p=
-(v \wedge \overline{ \imath \psi \wedge \ns w}, \Omega)_p=\]
\[=-(-1)^h(v\wedge \overline{ \ns \ns  (\imath \psi \wedge \ns w)},\Omega)_p=
-(-1)^h(v\wedge  \overline{ \ns (\ns \Psi  \ns )(w)},\Omega)_p \]
On the other hand,
\[(\Psi (v), w)_p=(v,\Psi^*(w))_p=(v\wedge \overline {\ns \Psi^*(w)}, \Omega)_p\]
This implies 
\[\ns \Psi^*=-(-1)^h\ns (\ns \Psi  \ns )\]
which is   equivalent to \[\Psi^*=-(-1)^h\ns \Psi  \ns \]
that is  the relation  we wanted to check.
The adjoint of this  equation immediately proves the relation also for the generators  \(\imath L_{\bk,\bj} \), (\(j\neq k\)) and   \(\imath A_i\).
\qed\\

\begin{teo}
\label{teo:ggru}
The Lie algebra bundle  \(\ggru\) is the full real Lie subalgebra bundle of $\gg$ of operators which preserve the form $<~,~>$, and its fibre is isomorphic to \(\mathbf{su}(s+1, s+1)\).
\end{teo}
\dimo
As usual, let us fix once and for all a point $p\in X$, which will be omitted from the notation when not strictly necessary.\\
In view of part $c)$ of Theorem \ref{teo:slstructure} we have to compute the signature of the form \(< \ , >_p\) when restricetd to \({\mathcal I}_{-s}\). It is convenient to use a basis \(\{c_r\}\) of  \({\mathcal I}_{-s}\)    which differs from the basis  \(\{b_r\}\) provided in the proof of Theorem  \ref{teo:slstructure} only for some signs.\\
Namely,   
  \(\{c_r\}\) (\(0\leq r\leq 2s +2\)) is given by the following monomials:
\begin{itemize}
\item  \(c_r= b_r = w_{0}\wedge \ldots \wedge \widehat {w_r} \wedge \cdots\wedge w_{s}\), for \(r \in \{0,1,2,\ldots , s\}\);
\item  \(c_{s+1+r}=w_{0}\wedge \cdots\wedge w_r\wedge  \overline{w}_r \wedge  \cdots \wedge w_{s}\), for \(r \in \{0,1,2,\ldots , s\}\).

\end{itemize}
By construction, for every \(r<j\),   \(<c_r, c_j>_p=0\) unless \(j=s+1+r\) and in this case we have that  \[<c_r, c_{s+1+r}>_p= \imath^{s(s+2)}(w_{0}\wedge \ldots \wedge \widehat {w_r} \wedge \cdots\wedge w_{s}\wedge \overline { w_{0}\wedge \cdots\wedge w_r\wedge  \overline{w}_r \wedge  \cdots \wedge w_{s}}, \Omega)_p= \] \[=\imath^{s(s+2)} (w_{0}\wedge \ldots \wedge \widehat {w_r} \wedge \cdots\wedge w_{s}\wedge \overline { w_{0}}\wedge \cdots\wedge \overline{w_r}  \wedge  w_r  \wedge  \cdots \wedge \overline { w_{s}}, \Omega)_p=\] 
\[=\imath^{s(s+2)}(-1)^{s+1} (w_{0}\wedge \ldots \wedge w_r \wedge \cdots\wedge w_{s}\wedge \overline { w_{0}}\wedge \cdots\wedge \overline{w_r}   \wedge  \cdots \wedge \overline { w_{s}}, \Omega)_p=\] 
\[=\imath^{s(s+2)}(-1)^{s+1} (-1)^s\cdots (-1)^0 (w_{0}\wedge\overline{w_0}\wedge  \ldots \wedge w_r \wedge \overline{w_r }\wedge \cdots\wedge w_{s}\wedge \overline { w_{s}}, \Omega)_p=\] 
\[=\imath^{s(s+2)}(-1)^{s+1} (-1)^s\cdots (-1)^0 (-\imath )^{s+1}=\imath^{s(s+2)}(-1)^s\cdots (-1)^0  \imath ^{s+1}=\] 
\[=\imath^{s(s+2)}(-1)^{\frac{s(s+1)}{2}} \imath ^{s+1}=\imath^{2s^2+1}\]  
Thus we notice that \(<c_r, c_{s+1+r}>_p\) does not depend on the index \(r\), being equal to $\imath$ when $s$ is even and to $-\imath$ when $s$ is odd.  If follows that  
 for every \(s\) the  signature   of \(< \ , \ >_p \)  is \((s+1,s+1)\).
 From  Theorem \ref{teo:antihermitean} and the above remark on the signature one deduces that the fibre of  \(\ggru\) at $p$ can be identified with a subalgebra  of \(\mathbf{su}(s+1,s+1)\subseteq \mathbf{sl}(2s +2,\C)\cong \gg_{,p}\). Since,  by construction, \(\ggru\) is a real form of \(\gg\),  we can replace   \(\subseteq\) with  \(=\) in the  inclusion above  and the claim follows.
\qed\\

This   theorem  also provides us the key ingredient to understand the structure of the bundle  $\ggr$, which is the most natural real Lie algebra bundle associated to the (almost, pointwise) $s$-\ka structure.  Indeed, it is generated at every point by the values of the operators $L_{ij}$, $V_k$ and their pointwise adjoints. The main tool will be  a new hermitean inner product $<<~,~>> $ defined starting from $<~,>$ on  the complex bundle of vector spaces $\bigwedge^*_\C T^*X$ which we now introduce:

\begin{dfn}
\label{dfn:secondhermitean}
For every $p\in X$ there is a natural  non degenerate Hermitean inner product $<<~,~>>_p$ on $\bigwedge_\C^*T^*_pX$, defined, on homogeneous elements \(\alpha, \beta\),  as:
\[ <<\alpha,\beta>>_p = \imath^{s+1+deg \ \beta }<\alpha, \beta>_p\]
We indicate with $<<~,~>>$ the corresponding form with values in smooth functions.
\end{dfn}

\begin{teo}
\label{teo:ggr}
The Lie algebra bundle  \(\ggr\) is the full real Lie subalgebra bundle of $\gg$ of operators which preserve the form $<<~,~>>$, and its fibre is isomorphic to \({\bf su}(s+1,s+1)\).
\end{teo}
\dimo
As usual, let us fix once and for all a point $p\in X$.
From  Proposition \ref{primeprop} we immediately deduce  that  the  weight components $\mathcal{I}_k$ are mutually $\quad <<~,~>>_p$-orthogonal and  $\quad <<~,~>>_p$ is nondegenerate when restricted to any one of them.

Therefore, for  the first part of the claim it suffices to show that the generators of \(\ggr\) preserve \(<< \ , \ >>\), then a dimensional argument can be applied (since \(\ggr\)  is a  proper real form of \(\gg\)).

Let \(\Gamma\) be any one of the   generators \(L_{jk},V_j,\Lambda_{jk},A_j\)  of \(\ggr\) . Then \(\imath \Gamma\) is a generator of \(\ggru\) and, by Theorem \ref{teo:ggru}, given two homogeneous elements  \(\alpha, \beta\in \bigwedge_\C^*T^*_pX\), it satisfies:
\[<\imath \Gamma (\alpha) , \beta>_p+ <\alpha, \imath \Gamma(\beta)>_p=0\]
Therefore,
\[<< \Gamma(\alpha) , \beta>>_p+ <<\alpha, \Gamma(\beta) >>_p=\] \[=\imath^{s+1+deg ~ \beta }<\Gamma (\alpha) , \beta>_p+\imath^{s+1+deg \  \beta+ deg \  \Gamma } <\alpha, \Gamma(\beta)>_p =\]
\[=\imath^{s+deg \ \beta }\left (  <\imath \Gamma (\alpha) , \beta>_p- \imath^{ deg \ \Gamma } <\alpha, \imath  \Gamma(\beta)>_p   \right ) =\]
\[=\imath^{s+deg \ \beta }\left (   < \imath \Gamma (\alpha) , \beta>_p+ <\alpha, \imath  \Gamma(\beta)>_p   \right ) =0\]
since \(deg \ \Gamma\) is equal to \(2\) or \(-2\).

For the second part of the claim, it suffices to compute the signature of the form \(<< \ , >>_p\) when restricetd to \({\mathcal I}_{-s}\).  Using  the  basis \(\{c_r\}\) of  \({\mathcal I}_{-s}\)    introduced in the proof of Theorem \ref{teo:ggru} we have:
\[<<c_r, c_{s+1+r}>>_p= \imath^{2s+3 }<c_r, c_{s+1+r}>_p=\imath^{2s^2+2s+4}=\imath^{2s(s+1)}=1\]
which shows that  the total signature is \((s+1,s+1)\).
\qed\\

We want now to characterize  explicitely the matrices of \(\ggr_{,p}\) with respect to the natural basis \(\{c_i\}\) of $\mathcal{I}_{-s}$ defined in the proof of Theorem \ref{teo:ggru}. Notice that the basis is not real, but we will show that the matrices are nevertheless in the standard form for $\mathbf{su}(s+1,s+1)$. 
\begin{pro}
\label{pro:matrici} With respect to the basis \(\{c_i\}\) of $\mathcal{I}_{-s}$, the algebra \(\ggr_{ , p}\)  is faithfully presented  as  the  algebra of matrices 
\[\left(
\begin{array}{cc}
  D & H_2  \\
H_1 &   -{\overline D}^t 
\end{array}
\right)\]
with   \(D\)  an $(s+1)\times (s+1)$ complex matrix  and \(H_1, H_2\)   two $(s+1)\times (s+1)$ complex antihermitean matrices.

\end{pro}

\dimo We start by noticing that  the operators of degree zero \(L_{j\bk}\) (\(j\neq k\)) (which lie in $\lefs_{p} \subset \ggr_{,p}$) can be expressed in terms of the basis of quadratic invariant monomials as:
\[2L_{j\bk}=[E_{w_j}, I_{w_k}]+ [E_{{\overline w}_j}, I_{{\overline w}_k}]\]

This allows us to compute, for \(i=0,1,2,\ldots , s\):
\[L_{j\bk}( c_i)=0 \qquad if \; i\neq j\]
\[L_{j\bk}(c_j)=\frac{1}{2}[E_{w_j}, I_{w_k}](c_j)=- I_{w_k}E_{w_j}(c_j)=- (-1)^{j+k}c_k\]
and 
\[L_{j\bk}( c_{s+1+i})=0 \qquad if \; i\neq k\]
\[L_{j\bk}(c_{s+1+k})=\frac{1}{2}[E_{{\overline w}_j}, I_{{\overline w}_k}](c_{s+1+k})=E_{{\overline w}_j} I_{{\overline w}_k}(c_{s+1+k})= (-1)^{j+k}c_{s+1+j}\]

This means that the matrices  of  the degree 0 subalgebra generated by the operators  \(L_{j\bk}\) have real coefficients and, more precisely  they are  all the matrices with the  following block-form:
\[
\left(
\begin{array}{cc}
  A  &  0   \\
  0   & -A^t   
\end{array}
\right)
\]
where \(A\) is a real \((s+1)\times (s+1)\) matrix with trace zero.  This explicitely establishes an isomorphism between \(<L_{j\bk}>_\R\) and \({\bf sl}(s+1, \R)\).\\
The computation of the matrices of  the operators  \(V_j\) is made easier by the use of the  relation contained in the proof of Proposition \ref{pro:ggrinvariant}:\[V_j=\frac{\imath}{2}  [E_ {w_j}, E_{{\overline w}_j}]\]
We can now observe that, for \(i=0,1,2,\ldots , s\):
\[V_j( c_i)=0 \qquad if \; i\neq j\]
\[V_j(c_j)=\frac{\imath }{2}2E_{w_j} E_{{\overline w}_j}(c_j)=\imath c_{s+1+j}\]
This, together with the observation that  
\[V_j( c_{s+1+i})=0 \qquad \forall i=0, \ldots , s\]
implies that the matrix of \(V_j\) has the following 
block-form:
\[
\left(
\begin{array}{cc}
  0  & 0   \\
  \imath B &   0   
\end{array}
\right)
\]
where \(B\) is a real and symmetric  (actually diagonal) \((s+1)\times (s+1)\) matrix.

It follows that all the  matrices of the above form are in \(\ggr_{,p}\) since they provide  an irreducible representation for the adjoint action of \(<L_{j\bk}>_\R \cong {\bf sl}(s+1, \R)\): notice that the action of a matrix with upper diagonal $A$ over one with lower left block $\imath B$ is as follows:
\[\imath B\to  -\imath(BA ~+~ ^tAB)\]

As for the operators   \(L_{jk}\) (\(j\neq k\)) of degree two, as it has been shown in Proposition \ref{pro:ggrinvariant}:
\[ 2L_{jk}=[E_{w_j},E_{\overline{w}_k}]-[E_{w_k},E_{\overline{w}_j}]\]
Therefore, 
\[L_{jk}( c_i)=0 \qquad if \; i\neq j,k \]
and 
\[L_{jk}(c_j)=-E_{{\overline w}_k}E_{w_j} (c_j)=(-1)^{j+k} c_{s+1+k}\]
\[L_{jk}(c_k)=E_{{\overline w}_j}E_{w_k} (c_k)=-(-1)^{j+k} c_{s+1+j}\]
This, together with the observation that  
\[L_{jk}( c_{s+1+i})=0 \qquad \forall i=0, \ldots , s\]
implies that the matrix of \(L_{jk}\) has the block-form:
\[
\left(
\begin{array}{cc}
  0  & 0   \\
  C &   0   
\end{array}
\right)
\]
where \(C\) is a real and antisymmetric  \((s+1)\times (s+1)\) matrix.
Then all the  matrices of the above form are in \(\ggr_{,p}\) since they provide  an irreducible representation for the action of \(<L_{j\bk}>_\R \cong {\bf sl}(s+1, \R)\)  similarly as before.

In the same way, acting with \(<L_{j\bk}>_\R \cong {\bf sl}(s+1, \R)\) on the adjoint operators \(\Lambda_{jk}\) and \(A_j\), we can show that \(\ggr_{,p}\) contains all the matrices of the form 
\[
\left(
\begin{array}{cc}
   0  & H  \\
  0 &   0   
\end{array}
\right)
\]
where  \(H\) is a complex antihermitean \((s+1)\times (s+1)\) matrix.

It is now immediate to check that the matrices constructed above generate the matrix algebra as in the claim. To conclude it is enough to use  point $c)$ of Theorem \ref{teo:slstructure}, which states that the restriction of $\ggrp$ to $\mathcal{I}_{-s}$ is faithful. 
\qed\\

The following corollary refines the result of Theorem \ref{pro:dn}:
\begin{cor}
\label{cor:orthogonal} 
For a fixed \(p\in X\), the orthogonal algebra  \(\lefs_p\cong \mathbf {so}(s+1,s+1)\) coincides with the intersection 
\( \ggr_{,p} \cap \ggrsp \).
\end{cor}
\dimo On one side, by Theorem \ref{pro:dn}, we  know that \(\lefs_p\) lies in \( \ggr_{,p} \cap \ggrsp\) and  is isomorphic to  \(\mathbf {so}(s+1,s+1)\). 

On the other side, let us consider the associated matrix algebras with respect to the basis \(\{c_i\}\) of $\mathcal{I}_{-s}$; the explicit presentation of \(\ggrsp\)  is made by real matrices (see the proof of Theorem \ref{pro:realform}), and the matrices of $\ggr_{,p}$ are as in Proposition  \ref{pro:matrici} above. Therefore:
\[dim_\R \  \left(\ggr_{,p} \cap\, \ggrsp\right) \leq dim_\R \ \{\hbox{subspace of} \ \ggr_{,p} \ \hbox{of  matrices with real coefficients}\} =  \]
\[ = (s+1)^2+ s(s+1)=dim_\R \  \mathbf {so}(s+1,s+1)\]

\section{Cotangent bundles, Tori and abelian varieties}
\label{sec:tori}
Let $M$ be a smooth Riemannian manifold with metric $\mathbf{h}$, and let 
\[X = \underbrace{T^*M\otimes_M\cdots\otimes_M T^*M}_{\mbox{$s$ times}}\]
We will show that $X$ is naturally almost $s$-\ka. First of all, clearly $\mathbf{h}$ induces naturally a Riemannian metric $\mathbf{g}$ on $X$. We then have to define the differential forms $\omega_{jk}$. These will come in two sets, with different constructions: the ones in which $j$ or $k$ is equal to zero and the other ones. The first ones are the simplest to define: is $\pi_j$ is the natural projection from $X$ to the $j$-th copy of $T^*M$, then we define 
\[\omega_{0j} = \pi_j^*\omega_{st}\]
where $\omega_{st}$ is the standard symplectic form on the cotangent bundle $T^*M$. A proof that with these forms $X$ becomes polysymplectic can be found in \cite{G1}, Example 2.3. The forms $\omega_{jk}$ when $j,k\not= 0$ will be defined by a different method. First, observe that using the Levi-Civita connection associated to the metric (induced by $\mathbf{g}$ on the cotangent bundle of $M$) we have a natural identification at any point $Q = (p,\phi_1,...,\phi_s)\in X$ 
\[T_{Q}X \cong T_p M\oplus \underbrace{T^*_p M\oplus\cdots\oplus T^*_pM}_{\mbox{$s$ times}}\]
Let us call $W_{jk}\subset T_{Q}X$ the direct sum of the $j^{th}$ and of the $k^{th}$ summands among the copies of $T^*_pM$ in the identification above:
\[T^*_pM\oplus T^*_pM\cong W_{jk}\subset T_QX\]
Using the metric, we can define $\omega_{ij}\in\bigwedge^2T^*_QX$ simply by defining a natural element in $\bigwedge^2 W_{jk}^*$. To do so, it is enough to observe that the identity (bundle) map from $T_pM$ to itself is an element  $Id \in T^*_pM\otimes T_pM$. 
This space is naturally isomorphic (using the metric $\mathbf{h}$) to $T_pM\otimes T_pM$ and this last space maps naturally to 
\[\bigwedge^2 \left(T_p M\oplus T_pM\right)\cong \bigwedge^2 W_{jk}^*\subset \bigwedge^2T^*_Q X\] 
where the last inclusion is again induced by the use of the metric. The proof that these forms satisfy the almost $s$-\ka condition is a simple direct computation. \\

The following example is a direct generalization of Example 2.7 of \cite{G2}. 
Let $\Gamma_0,...,\Gamma_s\subset\R^r$ be $s+1$ lattices, and let 
\[X = \R^r/\Gamma_0\times\cdots\times \R^r/\Gamma_s\]
Then $X$ has a natural structure of almost $s$-\ka manifold of rank $r$. Indeed, the (flat) metric is clear from the definition. We have also a natural choice of global coordinates 
\[\left\{y^j_i~|i\in\{1,...,r\}~,~j\in\{0,...,s\}\right\}\]
using which one can give directly the expressions for the forms:
\[\omega_{jk} = \sum_{i=1}^rdy^j_i\wedge dy^k_i\]
With the above definitions it is immediate to check that we have a almost $s$-\ka structure. More generally, one could take the above definitions of the $\omega_{jk}$ as forms on $(\R^r)^{ s+1}$ and define $X$ as the quotient of this almost $s$-\ka manifold by any (not necessarily maximal rank) lattice $\Gamma\subset (\R^r)^{ s+1}$:
\[X =  (\R^r)^{ s+1}/\Gamma\]
It is again immediate to check that in this way we obtain a almost $s$-\ka manifold of rank $r$, which is compact when $\Gamma$ is of maximal rank. As the lattice $\Gamma$ varies, we obtain different (and in general not isomorphic) almost $s$-\ka structures. As mentioned above, the metric being flat, these manifolds are actually a step higher in the rigidity ladder: they are $s$-\ka (see Example 8.3 of \cite{G1}).\\

The argument of Theorem 3.2 in \cite{G2} can be used to produce many more examples of almost $s$-\ka manifolds, compact or otherwise, by performing fibred products of special lagrangian fibrations. If the structure is $s$-\ka (see \cite{G1}, Definition 7.2) then it is automatically almost $s$-\ka; $s$-\ka manifolds are however  extremely rigid and difficult to construct, even more than Calabi-Yau ones, and therefore although geometrically interesting they are certainly not the correct way to try and build almost $s$-\ka ones.

\begin{pro}
\label{pro:connpiatta} On all the bundles of Lie algebras $\gg,\ggr,\ggrs,\ggru,\lefs$ there are natural flat connections induced by the almost $s$-\ka structure and compatible with inclusions. In the special case of $s$-\ka manifolds this follows from the observation that all the natural generators of these bundles of algebras are parallel tensors with respect to the Levi-Civita connection.
\end{pro}
\dimo 
Once fixed on an open set $U$ an orientation of the bundle $W_0$, one can choose over $U$ determinations of all the operators $V_j, A_k$. These, together with the restrictions to $U$ of the global sections $L_{jk}, \Lambda_{jk}$, generate over $\R$ a (finite dimensional) Lie algebra which we define to be the set of flat sections over $U$ for our connection on the bundle $\ggr$. Of course, the set of restriction of these sections to a given point $p\in X$ is exactly $\ggr_{,p}$.  A different choice for the orientation of $W_0$ would determine a choice for the $V_j$ (and hence $A_j$) which differs at most by a sign, thus the Lie algebra generated will not change, and we have a well defined set of sections over $U$. It is immediate to check that in this way one obtains a locally constant sheaf of sections for $\ggr$, and this determines a flat connection. The argument for the other natural bundles is the same.
\qed

We now show that we have a representation of the flat sections of the bundles of  Lie algebras
$\gg,\ggr,\ggrs,\ggru,\lefs$ on the cohomology of an $s$-\ka manifold, 
induced by the representation on the space of forms. This will be done showing that the 
Laplacian $\Delta_{d}$ commutes with the action of 
generators of these spaces of sections, as in  Theorem 10.1 on page  46 of \cite{G1}.
\begin{teo}
\label{teo:skaid}
Let $(X,\omega_{1},...,\omega_{s},\mathbf{g})$ be a compact  orientable
s-\ka manifold. Then we have that if $\phi\in\{L_{jk}\}\cup\{V_j\}$, and $d$ is the de Rham differential:\\
1) $[\phi,d]~=~0$\\
2)If we define 
$d^{c}~:=~[\phi,d^{*}]$,
we have that
$dd^{c}~+~d^{c}d~=~0$;\\
3) $[\phi,\Delta_{d}]~=~[\phi^*,\Delta_{d}]~=~0$,
where $\Delta_{d}$ is the $d$-Laplacian relative to the metric 
$\mathbf{g}$ and to the orientation.
\end{teo}
\dimo
We adapt the proof of Theorem 10.1  of \cite{G1}.\\
1) This equation follows immediately from the fact that the forms $\omega_{jk}$ and the volume forms $Vol(W_j)$ of the distribtions $W_j$ are covariant constant with respect to the Levi-Civita connection, and therefore closed.\\
2)
If we write down the expression for $d^{c}$ in standard $s$-\ka coordinates centered 
at a point $p\in X$, we see that no derivative of the metric appears. Therefore, when we write down the expression   for 
$dd^{c}~+~d^{c}d$, only the first derivatives of the 
metric are involved. We skip the details, as they are completely analogous to those of, for example, ~\cite[Pages 111-115]{GH}.\\
It follows, as in the classical case of \ka manifolds, that  to prove the equation it is enough to 
reduce to the case of a constant metric. 
When the metric is flat, however, the 
equation is easily seen to be equivalent (using $1)$) to 
$[\phi,\Delta_d]=0$, which with a flat metric follows immediately from 
the fact that the two-form corresponding to $\phi$ is constant in flat (orthonormal) 
coordinates.\\
3) The second equation is the adjoint of the first. The first one, 
once written down explicitely in terms of $d$ and $d^*$, follows 
immediately from points $1)-2)$.

\qed
\begin{cor}
\label{cor:cohomology}
Let $(X,\omega_{1},...,\omega_{s},\mathbf{g})$ be a compact orientable s-\ka 
manifold. Then there is a canonical representation of the Lie 
algebras of flat global sections of the bundles $\gg,\ggr,\ggrs,\ggru,\lefs$ on $H^*(X,\mathbf{C})$. 
\end{cor}

\begin{teo}
\label{teo:rivfinito}
For a compact rank two $s$-\ka manifold $X$ with a global determination of $J$ the following are equivalent:\\
1) X  is  Calabi-Yau.\\
2) There are a complex torus $\mathbb{T}$ with a tranlation invariant $s$-\ka structure and a holomorphic covering map $f:\mathbb{T}\to X$ compatible with the $s$-\ka structures of $\mathbb{T}$ and $X$. 
\end{teo}
\dimo
In the direction from 1) to 2), the main point in the proof is the observation that with a choice of $J$ at $p\in X$, the assignment of any one-form $v_{10}\in W_0\subset T^*_p X$ is enough to determine a complete adapted coframe. Indeed, let us call such a form $v_{10}$, and let us call $v_{20}$ its image under $J$: 
\[v_{20} = Jv_{10}\]
Using the natural identifications between the various $W_j$ induced by the structure forms and the metric, one obtains then corresponding forms $v_{11},v_{1,2},...,v_{1s},v_{2s}$ and it is immediate to check that they form an adapted coframe at the point $p$. Now, if we use these forms to build the corresponding $w_{0},...,w_s$, whose wedge product is a nonzero form of type $(s+1,0)$. A different choice of the initial $v_{10}$ could differ from the first one by a rotation of angle $\theta$:
\[v_{10}^\prime = e^{2\pi\theta J}v_{10}\]
This then is reflected in a modification as follows in the holomorphic volume form:
\[w_0^\prime\wedge\cdots\wedge w_s^\prime = 
e^{2\pi(s+1)\theta \imath}w_0\wedge\cdots\wedge w_s\]
Assuming now to have determined a global holomorphic volume form $\Omega$, we see that for any point there are $s+1$ choices of $v_{10}$ which produce the equation
\[w_0\wedge\cdots\wedge w_s = \Omega\]
at the point $p$. In other words, the set of all possible choices inside $T^*X$ forms a $(s+1)$-sheeted covering $\tilde{X}$ of $X$, over which there are $s+1$ global sections of the pull back of the covering itself. Such a covering space moreover inherits all the local geometric properties of $X$: it is a compact $s$-\ka manifold of rank two, with a global determination of $J$ and of a holomorphic volume form, and with a global smooth form $v_{10}$ determining the holomorphic volume form following the procedure described above. Correspondingly, there is a {\em global} determination of smooth forms $v_{10},...,v_{2s}$ determining an adapted coframe at all points $p\in X$. It is clear that from the covariance of the holomorphic volume form under parallel transport one obtains that any corresponding smooth determination of this adapted coframe will be also covariant constant with respect to Levi-Civita, and therefore a global determination will determine $s(s+1)$ covariant constant and orthonormal differential forms (and, dualizing, vector fields). This immediately shows that $X$ must be metrically flat, and being compact we have that it must me a torus. The global covariant orthonormal vector fields have therefore associated to them local coordinates which become global ones on the $\R^{2(s+1)}$ which is the universal covering   space of $\tilde{X}$. This (and the fact that the forms associated to these coordinates determine an adapted coframe) guarantee that the $s$-\ka structure is translation invariant with respect to the natural translation operation on $\tilde{X}$.\\ 
In the direction from 2) to 1), notice that you can always write $\mathbb{T}$ as $\C^{s+1}/\Gamma$ (as a complex manifold), with $\Gamma$ a lattice. Then observe that  $\C^{s+1}$ admits a translation invariant $s$-\ka structure with the complex structure inducing the $J$ operator; this $s$-\ka structure is therefore induced on $\mathbb{A}$. This argument works equally well when one has a covering of an abelian variety.\qed\\

\begin{rmk}
\label{rmk:abelianska}
The proof of the theorem shows in particular the following: if $A$ is a complex abelian variety of (complex) dimension $s+1$, then it is possible to put on it an $s$-\ka structure, with a translation invariant \ka metric, and with the complex structure of $A$ giving a determination of the operator $J$.
\end{rmk}

The preceding remark  shows that one can think of these $s$-\ka structures as "decorations" or "enrichments" of the underlying \ka structure. As such, one can use them to study the moduli problems for Abelian varieties by first studying the moduli problems of related $s$-\ka manifolds.\\

\section{Moduli of pointed elliptic curves}
\label{sec:moduli}
In this section we make a first attempt to put in contact the theory of $s$-\ka manifolds with that of moduli of pointed elliptic curves. We feel that there is (still partly hidden) much bigger interaction, which we intend to study in the future. Here, for $s\geq 2$, we construct a pointwise $(s-1)$-\ka structure  (of rank $2$) on the (open) moduli space $M_{1,s+1}$ with a fixed \ka structure.\\  
Let $E$ be an elliptic curve, with punctures $p_1,...,p_{s+1}$. The tangent space at the point $[E,p_1,...,p_{s+1}]$ of $M_{1,s+1}$ is  $H^1(E, T_E(-\sum_i p_i))$, and there is a short exact sequence of coherent sheaves
\[0\to T_E(-\sum_i p_i) \to T_E \to T_E/T_E(-\sum_i p_i)\to 0\]
from which one obtains the following exact sequence of complex vector spaces
\[0\to H^0(E,T_E)\to \bigoplus_i T_{p_i} E\to T_{[E,p_1,...,p_{s+1}]}M_{1,s+1} \to H^1(E,T_E)\to 0\]
The first and the last one of these vector spaces have dimension one, while the two intermediate ones have both dimension $s+1$.  In particular, there are canonical inclusions $T_{p_i} E\to T_{[E,p_1,...,p_{s+1}]}M_{1,s+1}$, and moreover one has that (with respect to these inclusions)
\[V_{E,p_1,...,p_{s+1}} := Im(\bigoplus_i T_{p_i} E)\subset  T_{[E,p_1,...,p_{s+1}]}M_{1,s+1}\]
has codimension one, with a permutation invariant syzygy among the images of a set of generators induced by a single translation invariant global holomorphic vector field on $E$. We indicate with $V$ the associated sub-bundle of $TM_{1,s+1}$. Using the metric on both sides, and taking the exterior power, we obtain a natural map
\[{\bigwedge^2}_\R\left(\bigoplus_i T^*_{p_i} E\right)\to {\bigwedge^2}_\R T^*_{[E,p_1,...,p_{s+1}]}M_{1,s+1}\]
As $E$ is an elliptic curve, for any pair $p_j,p_k$ of its points one has a natural complex linear and isometric identification (given by translation by $p_j - p_k$) from  $T_{p_k}E \to T_{p_j}E$, which can be equivalently seen as an element of 
$T^*_{p_k}E \otimes_\R T_{p_j}E$, and also (using the metric) as an element of $T^*_{p_j}E \otimes_\R T^*_{p_k}E$. Using the natural map 
\[T^*_{p_j}E \otimes_\R T^*_{p_k}E\to \bigwedge^2_\R\left(T^*_{p_j}E \oplus_\R T^*_{p_k}E\right)\] 
we have obtained a natural element 
\[\omega_{E,p_j,p_k}\in {\bigwedge^2}_\R \bigoplus_i T_{p_i} E\]
Varying the curve and the points, and using the map between exterior powers described before, we obtain a natural two-form
\[\omega_{jk}\in \Omega^2_\R(M_{1,s+1})\]
These two-forms are all well defined, and we want to restrict the set of forms 
\[\omega_{jk}~:~j,k\in\{1,...,s\}\]
(excluding the index $s+1$)  to the bundle $V\subset TM_{1,s+1}$. Using the metric of $M_{1,s+1}$ one can define a complementary (complex dimension one) subbundle $C\subset TM_{1,s+1}$, and induce on it a metric. On $V$ instead one induces a metric using the forms $\omega_{jk}$, and therefore in the end  one obtains a pointwise (s-1)-\ka structure (this is a pointwise verification) on $M_{1,s+1}$. As mentioned at the beginning of the section, this should be only part of the story: for instance, as a first step, it should be possible to go up to a full nondegenerate almost $s$-\ka structure with a little more effort.\\

\end{document}